\input amstex
\magnification=1200
\documentstyle{amsppt}
\vsize 18.8cm

\def\dist{\operatorname {\roman {dist}}}

\def\eps{\epsilon}
\def\om{\Omega}
\def\co{\colon}
\def\R{\Bbb R}
\def\C{\Bbb C}
\def\Z{\Bbb Z}
\def\L{\Bbb L}
\def\d{\,\roman{d}}

\def\supp{\operatorname{\roman{supp}}}
\def\diverg{\operatorname{\roman{div}}}

\def \R{\Bbb R}

\def\cl{\overline}

\def\N{\Bbb N}
\def\L{\Bbb L}
\def\H{\Bbb H}
\def\loc{\roman{loc}}
\def\i{\roman i}
\newcount\fornumber\newcount\artnumber\newcount\tnumber
\newcount\secnumber
\def\Ref#1{%
  \expandafter\ifx\csname mcw#1\endcsname \relax
    \warning{\string\Ref\string{\string#1\string}?}%
    \hbox{$???$}%
  \else \csname mcw#1\endcsname \fi}
\def\Refpage#1{%
  \expandafter\ifx\csname dw#1\endcsname \relax
    \warning{\string\Refpage\string{\string#1\string}?}%
    \hbox{$???$}%
  \else \csname dw#1\endcsname \fi}

\def\warning#1{\immediate\write16{
            -- warning -- #1}}
\def\CrossWord#1#2#3{%
  \def\x{}%
  \def\y{#2}%
  \ifx \x\y \def\z{#3}\else
            \def\z{#2}\fi
  \expandafter\edef\csname mcw#1\endcsname{\z}
\expandafter\edef\csname dw#1\endcsname{#3}}
\def\Tag#1#2{\begingroup
  \edef\mmhead{\string\CrossWord{#1}{#2}}%
  \def\writeref{\write\refout}%
  \expandafter \expandafter \expandafter
  \writeref\expandafter{\mmhead{\the\pageno}}%
\endgroup}  
  
\openin15=\jobname.ref
\ifeof15 \immediate\write16{No file \jobname.ref}%
\else      \input \jobname.ref
\fi
\closein15
\newwrite\refout
\openout\refout=\jobname.ref
  
\def\newsection{\global\advance\secnumber by 1\fornumber=0\tnumber=0}

\def\dff#1..{%
     \global\advance\fornumber by 1
     \Tag{#1}{(\the\secnumber.\the\fornumber)}
 \unskip          \the\secnumber.\the\fornumber
   \ignorespaces
}%

\def\rff#1..{\Ref{#1}}
\def\dft#1..{%
     \global\advance\tnumber by 1
     \Tag{#1}{\the\secnumber.\the\tnumber}
 \unskip          \the\secnumber.\the\tnumber
   \ignorespaces
}%
\def\rft#1..{\Ref{#1}}

\def\dfa#1..{%
     \global\advance\artnumber by 1
     \Tag{#1}{\the\artnumber}
 \unskip          \the\artnumber
   \ignorespaces
}%
\def\rfa#1..{\Ref{#1}}
   
\def\dep#1..{\Tag{#1}{}}
\def\rep#1..{p. \Refpage{#1}}
     
\topmatter
\title
 Characterization of the limit
of some \\ higher dimensional thin domain problems\endtitle
\author Thomas Elsken --
Martino Prizzi\endauthor
\leftheadtext{T. Elsken and M. Prizzi}
\rightheadtext{Higher dimensional thin domain problems}
\address Universit\"at Rostock, 
Fachbereich Mathematik, Universit\"atsplatz 1, 18055 Rostock, Germany
\endaddress\email
thomas.elsken\@mathematik.uni-rostock.de\endemail
\address 
Universit\`a degli Studi di Trieste,
Dipartimento di Scienze Matematiche, Via Valerio 12/b, 34100 Trieste,
Italy
\endaddress \email
prizzi\@mathsun1.univ.trieste.it\endemail

\abstract
A reaction-diffusion equation on a family of three dimensional thin domains, collapsing 
onto a two dimensional subspace, is considered.
In \cite{\rfa pr..} it was proved that, as the thickness of the domains tends to zero, 
the
solutions of the equations converge in a strong sense to the solutions of an abstract 
semilinear parabolic equation living in a closed subspace of $H^1$. Also, existence and
upper semicontinuity of the attractors was proved. In this work, for a specific class of
domains, the limit problem is completely characterized as a system of two-dimensional
reaction-diffusion equations, coupled by mean of compatibility and balance boundary
conditions. 
\endabstract

\endtopmatter

\document\newsection\head \the\secnumber.
Introduction\endhead

Let $\Omega\subset\R^{N+M}$ be an open bounded domain with Lipschitz boundary.
Write $(x,y)$ for a generic point of $\R^{N+M}$.
For $\epsilon>0$, let us consider the `squeezing operator'
$T_\epsilon\colon\R^{N+M}\to\R^{N+M}$, $(x,y)\mapsto(x,\epsilon y)$, and define
$\Omega_\epsilon:=T_\epsilon(\Omega)$. Let $\Gamma$ be a relatively closed
portion of $\partial\Omega$ and let $\Gamma_\epsilon:=T_\epsilon(\Gamma)$.
Let us consider the following reaction-diffusion equation
$$ \cases
u_t=\Delta u+
f(u), &t>0,\, 
 (x,y)\in\Omega_\eps
\\ \partial _{\nu_\eps} u=0,
& t>0,\, (x,y)\in\partial\Omega_\eps\setminus\Gamma_\eps
\\ u=0,
& t>0,\, (x,y)\in\Gamma_\eps.
\endcases
\tag \dff e1.. 
$$
Here $\nu_\eps$ is the exterior normal vector field on $\partial \om_\eps$.
We assume that $f$ satisfies the following condition:\medskip
(H1) $f\in C^1(\R\to\R)$ and
$|f'(s)|\le C(|s|^\beta+1)$ for $s\in \R$, where $C$ and
$\beta\in\left[0,\infty\right[$ are arbitrary real constants. If
$n:=M+N>2$ then  in addition, $\beta\le (p^*/2)-1$, where
$p^*=2n/(n-2)>2$.\medskip 

Let $H^1_{\Gamma_\epsilon}(\Omega_\epsilon)$ be the closure in
$H^1(\Omega_\epsilon)$ of the space of all $C^1(\cl\om_\eps)$-functions
such that $u=0$ on $\Gamma_\eps$.
Then it is well known that  
equation \rff e1.. generates a 
semiflow $\tilde\pi_\eps$ on $H^1_{\Gamma_\eps}(\om_\eps)$.
If we suppose in addition that $f$ satisfies the dissipativeness 
condition\medskip
(H2) $ \limsup_{|s|\to \infty}f(s)/s\le -\zeta$ for some
$\zeta>0$,\medskip
then the semiflow $\tilde\pi_\epsilon$ is defined for all $t\geq0$ and it 
posseses a compact global attractor $\tilde{\Cal A}_\eps$.

As $\eps\to 0$ the thin domain $\om_\eps$ degenerates to an
$N$-dimensional domain.
Then the question arises, what happens in the limit to the family
$(\tilde\pi_\eps)_{\eps>0}$ of semiflows and to the family
$(\tilde\Cal A_\eps)_{\eps>0}$ of attractors. Does there exist a limit
semiflow and a corresponding limit attractor?

 This problem was first considered by Hale and
Raugel in \cite{\rfa HaRau1..} for the case when $M=1$ and the domain $\om$
is the {\it ordinate set\/} of a smooth positive function $g$
defined on an $N$-dimensional domain $\omega$, i.e. $$\om=\{\, (x,y)\mid
\text{$x\in\omega$ and $0<y<g(x)$}\,\},$$ with $\Gamma=\emptyset$
(resp. $\Gamma=\{\, (x,y)\mid
\text{$x\in\partial\omega$ and $0<y<g(x)$}\,\}$).

The authors prove that, in
this case, there exists a limit semiflow $\tilde\pi_0$,  which is
defined by the $N$-dimensional boundary value problem 
$$
\cases u_t=(1/g)\diverg(g\nabla u)+ f(u), &t>0,\,
 x\in\omega
\\{{\partial u}\over{\partial\nu}} u=0\quad\text{(resp. $u=0$}), & t>0,\, 
x\in\partial\omega. \endcases\tag \dff 
e2..
$$ 
Moreover, $\tilde\pi_0$ has a global attractor $\tilde \Cal A_0$ and, in
some sense, the family $(\tilde \Cal A_\eps)_{\eps\ge 0}$ is
upper-semicontinuous at $\eps=0$. See also \cite{\rfa Rau1..} and the
rich bibliography contained therein.

 If
the domain $\om$ is not the ordinate set of some function (e.g. if
$\om$ has holes or different horizontal branches) then \rff e2.. can
no longer be a limiting equation for \rff e1... Nevertheless, K. Rybakowski and the second 
author proved in \cite{\rfa pr..} that the family $\tilde \pi_\eps$
still has a limit semiflow. Moreover, there exists a limit global
attractor and the upper-semicontinuity result continues to hold.

In order to describe the main results of \cite{\rfa pr..}
 we first transfer the family \rff e1.. to boundary value problems on the fixed
domain $\om$. More explicitly, we use the linear isomorphism
$\Phi_\eps\co H^1(\om_\eps)\to H^1(\om)$,	$u\mapsto u\circ
T_\eps$, to transform problem \rff e1.. to the equivalent problem
$$ \cases u_t=\Delta_xu+\frac 1{\eps^2}\Delta_yu+ f(u),
&t>0,\,
 (x,y)\in\Omega\\
 \nabla_x u\cdot\nu_x+\frac 1{\eps^2}\nabla_y u\cdot\nu_y=0, & t>0,\,
(x,y)\in\partial\Omega\setminus\Gamma\\	 
u=0,&  t>0,\,
(x,y)\in\Gamma\\
\endcases\tag \dff e3.. $$ on  $\om$.
Here, $\nu=(\nu_x, \nu_y)$ is the exterior normal vector field on
$\partial \om$.

Let $H^1_{\Gamma}(\Omega)$ be the closure in
$H^1(\Omega)$ of the space of all $C^1(\cl\om)$-functions
such that $u=0$ on $\Gamma$.
Then equation \rff e3.. can be written in the abstract form
$$\dot u+A_\eps u=\hat f(u)\tag\dff e4..$$ where $\hat f\co H^1_\Gamma(\om)\to
L^2(\om)$ is the Nemitski operator
 generated by the function $f$, and $A_\eps$
is the selfadjoint linear operator (with compact resolvent) induced by the
following bilinear form 
$$a_\eps (u,v):=\int_\om(\nabla_xu\cdot\nabla_xv+\frac
1{\eps^2}\nabla_yu\cdot\nabla_yv)\,\d x\,\d y, 
\quad u, v\in H^1_\Gamma(\om).$$

Equation \rff e4.. then defines a semiflow $\pi_\eps$ on $H^1_\Gamma(\om)$
which is equivalent to $\tilde \pi_\eps$ and has the global
attractor $\Cal A_\eps:=\Phi_\eps(\tilde\Cal A_\eps)$, consisting
of the orbits of all full bounded solutions of \rff e4...

 Notice
that, for	every fixed $\eps>0$ and $u\in H^1_\Gamma(\om)$, the formula
$$|u|_\eps=\left(a_\eps(u,u)+|u|^2_{L^2(\om)}\right)^{1/2}$$
 defines a norm on $H^1_\Gamma(\om)$ which is equivalent to 
 $|\cdot|_{H^1_\Gamma(\om)}$.
However, $|u|_\eps\to \infty$ as $\eps\to 0^+$ whenever
$\nabla_yu\not=0$ in $L^2(\om)$.

In fact, we see that for $u\in H^1_\Gamma(\om)$ $$\lim_{\eps\to
0^+}a_\eps(u,u)=\cases \int_\om |\nabla_xu|^2\,\d x\,\d y,&\text{if
$\nabla_yu=0$}\\+\infty,& \text{otherwise.}\endcases$$ Thus the family
$a_\eps(u,u)$, $\eps>0$, of real numbers has a finite limit (as
$\eps\to 0$) if and only if $u\in H^1_{\Gamma,s}(\om)$, where we define
$$H^1_{\Gamma,s}(\om):=\{\,u\in H^1_\Gamma(\om)\mid \nabla_yu=0\,\}.$$ 
This is a closed
linear subspace of $H^1_\Gamma(\om)$. 

The corresponding limit bilinear form is given by the formula:
$$a_0 (u,v):=\int_\om \nabla_xu\cdot\nabla_xv\d x\,\d y, \quad u, v\in
H^1_{\Gamma,s}(\om).\tag\dff lbil..$$

Assume from now on that 
$H^1_{\Gamma,s}(\om)$ is infinite dimensional.
Then the form $a_0$ uniquely determines a densely defined selfadjoint
linear operator $$A_0\co D(A_0)\subset H^1_{\Gamma,s(\om)}\to 
L^2_{\Gamma,s}(\om)$$
by the usual formula $$a_0(u,v)=\langle
A_0u,v\rangle_{L^2(\om)},\quad \text{ for $u\in D(A_0)$ and $v\in
H^1_{\Gamma,s(\Omega)}$.}\tag\dff lop..$$ Notice that $A_0$ has compact
resolvent. 
Here, $L^2_{\Gamma,s}(\om)$ is the closure of 
$H^1_{\Gamma,s}(\om)$ in the
$L^2$-norm, so $L^2_{\Gamma,s}(\om)$ is a closed linear subspace of
$L^2(\om)$.

One can show that the Nemitski operator $\hat f$ maps the space
$H^1_{\Gamma,s}(\om)$ into $L^2_{\Gamma,s}(\om)$. Consequently the abstract
parabolic equation 
$$
\dot u+A_0 u=\hat f(u)
\tag\dff e5..
$$ 
defines a
semiflow $\pi_0$ on the space $H^1_{\Gamma,s}(\om)$. This is the limit
semiflow of the family $\pi_\eps$. The following results
are proved in \cite{\rfa pr..}:
\proclaim{Theorem A} Let
$(\eps_n)_{n\in\N}$ be an arbitrary sequence of positive numbers
convergent to zero and $(u_n)_{n\in\N}$ be a sequence in
$L^2(\om)$ converging in the norm of $L^2(\om)$ to some $u_0\in
L^2_{\Gamma,s}(\om)$. Moreover, let $(t_n)_{n\in\N}$ be an arbitrary
sequence of positive numbers converging to some positive number
$t_0$.

Then
$$\bigl|e^{-t_nA_{\eps_n}}u_n-e^{-t_0A_{0}}u_0\bigr|_{\eps_n}\to
0\quad \text{as $n\to \infty$.}$$ If, in addition,	$u_n\in H^1(\om)$
for every $n\in \N$ and if
 $u_0\in H^1_{\Gamma,s}(\om)$, then
$$|u_n\pi_{\eps _n} t_n-u_0\pi_0 t_0|_{\eps_n}\to 0\quad \text{as
$n\to \infty$.}$$\endproclaim
 The limit semiflow $\pi_0$ possesses a global attractor
$\Cal A_0$. The upper-semicontinuity result alluded to above reads
as follows: \proclaim{Theorem B} The family of attractors
$\left(\Cal A_\eps\right)_{\eps\in[0,1]}$ is upper-semicontinuous
at $\eps=0$ with respect to the family of norms $|\cdot|_\eps$.

This means that $$\lim_{\eps\to 0^+}\sup_{u\in \Cal
A_\eps}\inf_{v\in \Cal A_0}|u-v|_\eps=0.$$
In particular, there exists an $\eps_1>0$
 and an open bounded set $U$ in $H^1(\om)$ including all the
attractors $\Cal A_\eps$, $\eps\in [0,\eps_1]$. \endproclaim

\remark{Remark}Theorems A and B were actually proved in the case $\Gamma=\emptyset$,
but the proof is valid (with only minor changes) also in the general case,
as long as $H^1_{\Gamma,s}(\om)$ is infinite dimensional.\endremark

The definition of the linear operator $A_0$, as given above, is
not very explicit. If $N=M=1$, however,  it was shown in \cite{\rfa pr..}
and \cite{\rfa pr2..} that
there is a large class of the so-called {\it nicely decomposed\/}
 domains on which $A_0$ can be characterized as a
system of one-dimensional second order linear differential operators,
coupled to each other by certain compatibility and Kirchhoff type
balance conditions. In this case, the abstract limit  equation
\rff e5.. is equivalent to a parabolic equation on a finite graph.
Roughly speaking, a planar domain $\om$ admits a {\it nice decomposition} if, 
up to a set of measure zero
contained in a set $Z$ of finitely many vertical lines, $\om$ can be
decomposed into finitely many domains $\om_k$, $k=1$, \dots,
$r$ in such a way that at $Z$ the various sets $\om_k$ and $\om_l$
`join' in a nice way. Points of $\cl\om\cap Z$ are, intuitively speaking,
those at which	connected components of the vertical sections
$\om_x$ bifurcate (see Figure 3 in \cite{\rff pr..}).
In higher dimensions it is not clear wether it is possible to describe 
a reasonable, sufficiently large, class of domains for which an explicit
characterization of $H^1_{\Gamma,s}(\om)$ and of $D(A_0)$ can be carried on.
Nevertheless, in some concrete cases, one can go along the same 
ideas of \cite{\rfa pr..} and give a nice characterization of
these spaces. In this paper we concentrate on 
the case $N=2$ and $M=1$ and we illustrate with two examples how this is 
possible. Our examples deal with a set $\om$ which is obtained by removing from
a cylinder a smaller cylinder contained in the interior of the first.
More precisely, take open sets $\omega$, $\omega_1$, $\omega_2$ and $\omega_3$ in $\R^2$ such
that 
$$
\alignedat1
&\omega \quad\text{is bounded, connected and has $C^2$ boundary,}\\
&\omega_2=\omega_3\subset\subset\omega\quad\text{have $C^2$ boundary,}\\
&\omega_1:=\omega\setminus \cl\omega_2.
\endalignedat
$$
Notice that $\omega_1$ is not necessarily connected. Moreover, let $h_1$, $h_2$ and $h_3$ be
positive real numbers, with
$h_1>h_2+h_3$. Then we define 
$$
\om:=\left(\omega\times]0,h_1[\right)\setminus
\cl{\left(\omega_2\times]h_3,h_1-h_2[\right)}.\tag\dff cili..
$$
Figure 1 below represents the domain $\Omega$, when $\omega$ and $\omega_2$ are balls
centered at $0$.

\hskip2.5cm\vbox to 3in{\vss\includegraphics{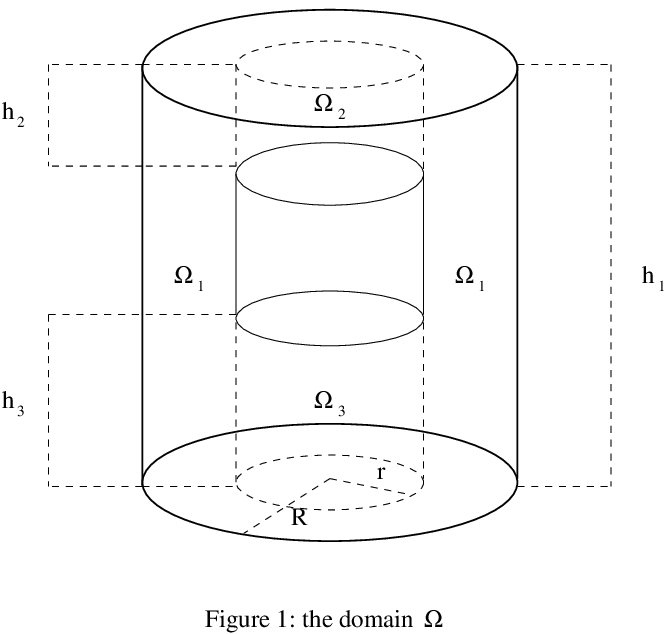}}

\bigskip

For later use we need also to define
$$
\alignedat1
\om_1:=&\omega_1\times]0,h_1[,\\
\om_2:=&\omega_2\times]h_1-h_2,h_1[,\\
\om_3:=&\omega_3\times]0,h_3[\\
\endalignedat
$$
and
$$
\alignedat1
\om_4:=&\omega_1\times]h_3,h_1-h_2[,\\
\om_5:=&\omega\times]h_1-h_2,h_1[,\\
\om_6:=&\omega\times]0,h_3[,\\
\om_7:=&\R^2\times]h_3,h_1-h_2[.\\
\endalignedat
$$
Finally, we set
$$
\Gamma_1:=\partial\omega\times[0,h_1], \quad
\Gamma_2:=\partial\omega_2\times[h_3,h_1-h_2],\quad
\Gamma_L:=\Gamma_1\cup\Gamma_2.
$$
We shall consider equation \rff e1.. on $\om_\eps=T_\eps(\om)$, where
$\om$ is the domain defined above,
with two different sets of boundary conditions, namely with $\Gamma=\emptyset$
and with $\Gamma=\Gamma_L$. We shall see that these different boundary 
conditions give rise to completely different behaviors as $\eps\to 0$.
In fact, when $\Gamma=\emptyset$, i.e. we impose the Neumann boundary condition
on the whole $\partial\om_\eps$, equation \rff e5.. is equivalent to the 
following system of two-dimensional reaction-diffusion equations
$$ \cases
u_{it}=\Delta u_i+
f(u_i), &t>0,\, 
 x\in\omega_i,\quad i=1,2,3,\\
u_1(x)=u_2(x)=u_3(x),&t>0,\,x\in\partial\omega_2,\\  
 \partial _{\nu_1} u_1=0,&
 t>0,\, x\in\partial\omega,\\ 
\sum_{i=1}^3h_i\nabla u_i\cdot\nu_i=0,
& t>0,\, x\in\partial\omega_2.
\endcases
\tag \dff e6.. 
$$
Here $\nu_i$, $i=1,2,3$, is the outward normal vector field on 
$\partial\omega_i$ for $i=1,2,3$ respectively. Observe that the three 
equations in \rff e6.. are coupled by compatibility and Kirchoff type balance
conditions on the `interface' $\partial\omega_2$. Figure 2 below illustrates the `limit'
of the family $(\Omega_\epsilon)$ as $\epsilon\to 0$ for the domain represented
in Figure 1.

\hskip1.5cm\vbox to 2in{\vss\includegraphics{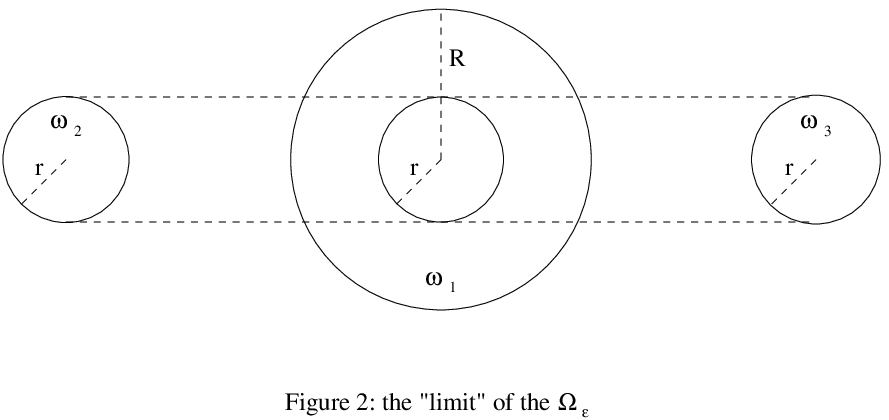}}
 
\bigskip

On the other hand, when $\Gamma=\Gamma_L$, i.e. we impose the Dirichlet
boundary condition on the `lateral' surface $\Gamma_L$, equation
\rff e5.. is equivalent to the 
following system of two-dimensional reaction-diffusion equations
$$ \cases
u_{it}=\Delta u_i+
f(u_i), &t>0,\, 
 x\in\omega_i,\quad i=1,2,3,\\
u_i(x)=0,&t>0,\,x\in\partial\omega_i,\quad i=1,2,3.\\  
\endcases
\tag \dff e7.. 
$$
So in this case the `limit' problem is a completely decoupled system
of scalar reaction-diffusion equations.

These two examples furnish a prototype for many concrete situations
that may occur in practice. 
In particular, we point out that the core of the problem consists in
proving regularity of the solutions of the linear equation
$$
A_0u=w,\quad\text{with $w\in L^2_{\Gamma,s}(\om)$.}
$$
Once the spaces  $L^2_{\Gamma,s}(\om)$, $H^1_{\Gamma,s}(\om)$
and $D(A_0)$  have been characterized, 
one can easily show that \rff e5.. is equivalent to a system of concrete
reaction-diffusion equations of type \rff e6.. or \rff e7...

Finally, as we shall explain in Section 3, the characterization of
$A_0$ and of its domain can be exploited to compute the eigenvalues of
$A_0$ in some specific situations, like the one represented in Figure 1. 
Of course, informations on the location and on the multiplicity of the eigenvalues of 
$A_0$ are very important in the study of local bifurcations of \rff e5...

\newsection\head \the\secnumber.
Characterization of $H^1_{\Gamma,s}(\om)$\endhead
We begin by recalling a general notion introduced in \cite{\rfa pr..}: we
say
that an open set
$\om\in
\R^{N+M}$ {\it has connected vertical sections\/} if for every $x\in \R^N$
the
$x$-section
$\om_x$ is connected. Of course, this section is nonempty if and only if
$x\in
P(\om)$, where
$P\co
\R^N\times \R^M\to
\R^N$,
$(x,y)\mapsto x$ is the projection onto the first $N$ components.
The following proposition was proved in \cite{\rfa pr..}:
\proclaim{Proposition \dft co..} Suppose $\om$ has connected vertical
sections.
 Let
$J:=P(\om)$ and
define the function $p\co J\to \left]0,\infty\right[$ by $x\mapsto
\mu_M(\om_x)$. If
  $u\in L^2(\om)$ satisfies $\nabla_y u=0$ in the distributional
sense, then there is a null set $S$ in $\R^{N+M}$ and a function $v\in
L^1_{\loc}(J)$ such that
$u(x,y)=v(x)$ for every $(x,y)\in \om\setminus S$. Moreover,
$p^{1/2}v\in L^2(J)$. If $u\in H^1(\om)$ then
$\partial_{x_i}v \in L^1_{\loc}(J)$ for $i=1$, \dots, $N$ and
 we can choose the null set $S$ so
that $u(x,y)=v(x)$ and $\partial_{x_i}u(x,y)=\partial_{x_i}v(x)$
for every $i=1$, \dots, $N$ and $(x,y)\in \om\setminus S$. Moreover,
$p^{1/2}\partial_{x_i}v\in L^2(J)$ for every $i=1$, \dots, $N$.
\qed\endproclaim
Now we come back to the domain $\om$ defined by \rff cili...
In what follows, we may assume indifferently that $\Gamma=\Gamma_L$
or $\Gamma=\emptyset$. For $k=1$, \dots, $7$ let us define
$$
H^1_s(\om_k):=\{\,u\in H^1(\om_k)\mid\nabla_y u=0~\text{a.e}.\,\}.
$$
Moreover, let us define $L^2_s(\om_k)$ as the closure of $H^1_s(\om_k)$
in $L^2(\om_k)$.
\proclaim{Lemma \dft restr..}
For $k=1$, \dots, $6$,
the following properties hold:
\roster
\item whenever $u\in L^2_{\Gamma,s}(\Omega)$, then $u|_{\Omega_k}\in
L^2_s(\Omega_k)$;
\item whenever $u\in H^1_{\Gamma,s}(\Omega)$, then $u|_{\Omega_k}\in
H^1_s(\Omega_k)$.
\endroster
\endproclaim
\demo{Proof}
Part (2) is obvious and part (1) follows directly from part (2) and from
the definition of $L^2_{\Gamma,s}(\Omega)$ and $L^2_s(\Omega_k)$.
\qed\enddemo

For $k=1,2,3$, let us define the spaces
$$
L_k:=L^2(\omega_k)\quad\text{and}\quad H_k:=H^1(\omega_k).
$$
Define on $L_k$ and $H_k$ the scalar products
$$
\langle u,v\rangle_{L_k}:=\int_{\omega_k}h_ku(x)v(x)\d x
$$
and
$$
\langle u,v\rangle_{H_k}:=\int_{\omega_k}h_ku(x)v(x)\d x
+\int_{\omega_k}h_k\nabla u(x)\cdot\nabla v(x)\d x
$$
respectively. Moreover, for $k=1,2,3$,  let us define the mapping
$$
\imath_k\co L^2_s(\om_k)\to L_k,\quad u\mapsto v,
$$
where $v$ is the function given by proposition \rft co... It turns out
that $\imath_k$ is an isometry of $L^2_s(\om_k)$ onto $L_k$ for
$k=1,2,3$.
Moreover, $\imath_k$ restricts to an isometry of $H^1_s(\om_k)$ onto
$H_k$
for $k=1,2,3$.
Let us define the product spaces
$$
L_\oplus:=L_1\oplus L_2\oplus L_3:=\{\,[u]=(u_1,u_2,u_3)\mid u_k\in
L_k,k=1,2,3\,\}
$$
and
$$
H_\oplus:=H_1\oplus H_2\oplus H_3:=\{\,[u]=(u_1,u_2,u_3)\mid u_k\in
H_k,k=1,2,3\,\}
$$
with the scalar products
$$
\langle [u],[v]\rangle_{L_\oplus}:=\langle u_1,v_1\rangle_{L_1}+
\langle u_2,v_2\rangle_{L_2}+\langle u_3,v_3\rangle_{L_3}
$$
and
$$
\langle [u],[v]\rangle_{H_\oplus}:=\langle u_1,v_1\rangle_{H_1}+
\langle u_2,v_2\rangle_{H_2}+\langle u_3,v_3\rangle_{H_3}
$$
respectively. It is easy to check that $L_\oplus$ and $H_\oplus$ 
are Hilbert spaces. Besides, let us define the map
$$
\imath_\oplus\co L^2_{\Gamma,s}(\om)\to L_\oplus,
\quad \imath_\oplus
u:=(\imath_1(u|_{\om_1}),\imath_2(u|_{\om_2}),\imath_3(u|_{\om_3})).
$$
Observe that 
$$
\langle u,v\rangle_{L^2(\om)}=\langle\imath_\oplus u,\imath_\oplus
v\rangle_{L_\oplus}
\quad\text{for $u$ and $v\in L^2_{\Gamma,s}(\om)$}
$$
and
$$
\langle u,v\rangle_{L^2(\om)}+a_0(u,v)=\langle\imath_\oplus
u,\imath_\oplus v\rangle_{H_\oplus}
\quad\text{for $u$ and $v\in H^1_{\Gamma,s}(\om)$.}
$$
It follows by Lemma \rft restr.. that  $\imath_\oplus$ is an isometry of
of $L^2_{\Gamma,s}(\om)$ into $L_\oplus$ and that $\imath_\oplus$
restricts to
an isometry of $H^1_{\Gamma,s}(\om)$ into $H_\oplus$.
Finally, let us define
$$
\gather
H^0_\oplus:=
\left\{ [u]\in H_\oplus\mid u_k\in H^1_0(\omega_k)~\text{for
$k=1,2,3$}\right\}
\endgather
$$
and
$$
H^C_\oplus:=
\left\{ [u]\in H_\oplus\mid \null^\tau
u_1(x)=\null^\tau u_2(x)=\null^\tau u_3(x)~\text{${\Cal H}^1$-a.e.
on $\partial\omega_2$}\right\}
$$
where 
${\Cal H}^1$ is the one-dimensional Hausdorff measure in $\R^2$ and
$\null^\tau u_k$ is the trace of $u_k$ on $\partial\omega_k$ for
$k=1,2,3$. We call
$$
\null^\tau u_1(x)=\null^\tau u_2(x)=\null^\tau u_3(x)
\quad\text{${\Cal H}^1$-a.e. on $\partial\omega_2$}
\tag\dff compcond..
$$
the {\it compatibility condition} on $\partial\omega_2$.

Now we are able to characterize the spaces $H^1_{\Gamma,s}(\om)$ and
$L^2_{\Gamma,s}(\om)$:

\proclaim{Theorem \dft char1..}
The following properties hold:
\roster
\item $\imath_\oplus(L^2_{\Gamma,s}(\om))=L_\oplus$;
\item $\imath_\oplus(H^1_{\Gamma,s}(\om))=H^C_\oplus$ if
$\Gamma=\emptyset$ and $\imath_\oplus(H^1_{\Gamma,s}(\om))=H^0_\oplus$
if $\Gamma=\Gamma_L$.
\endroster
\endproclaim
\demo{Proof}
We begin by proving (2). Let $\Gamma=\Gamma_L$ or $\Gamma=\emptyset$
and let $u\in H^1_{\Gamma,s}(\om)$. Let $\imath_\oplus
u:=[v]=(v_1,v_2,v_3)$. We shall prove that 
$$\null^\tau
v_1(x)=\null^\tau v_2(x)=\null^\tau v_3(x)\quad\text{${\Cal H}^1$-a.e.
on $\partial\omega_2$.}\tag\dff compat.. 
$$
By the definition of $\imath_\oplus$ and by Proposition \rft co..,
there exists a null set $S\subset\R^3$ such that
$$ \alignedat2
u(x,y)&=v_k(x)&\quad&\text{for all $(x,y)\in\Omega_k\setminus
S$ and for $k=1,2,3$.}\\  \endalignedat 
$$
On the other hand, again by Proposition \rft co.., we can find
two functions $v_5$ and $v_6\in H^1(\omega)$ and we can choose the set $S$
in such a way that
$$ \alignedat2
u(x,y)&=v_l(x)&\quad&\text{for all $(x,y)\in\Omega_l\setminus
S$ and for $l=5,6$.}\\  \endalignedat 
$$
It follows that
$$
v_1(x)=v_5(x)=v_6(x)\quad\text{a.e. in $\omega_1$},
$$
$$
v_2(x)=v_5(x)\quad\text{a.e. in $\omega_2$}
$$
and
$$
v_3(x)=v_6(x)\quad\text{a.e. in $\omega_3$}.
$$
Define the functions $\tilde v_5$ and $\tilde v_6\co\omega\to\R$ by
$$
\tilde v_5(x):=\cases
v_1(x)&\text{if $x\in\omega_1$}\\v_2(x)&\text{if $x\in\omega_2$}\\
0&\text{otherwise}
\endcases
$$
and
$$
\tilde v_6(x):=\cases
v_1(x)&\text{if $x\in\omega_1$}\\v_3(x)&\text{if $x\in\omega_3$}\\
0&\text{otherwise}
\endcases
$$
It follows that $\tilde v_5=v_5$ and $\tilde v_6=v_6$ almost everywhere in
$\omega$ and hence $\tilde v_5$ and $\tilde v_6\in H^1(\omega)$. This in
turns implies \rff compat.. (see \cite{\rfa alt.., Lemma A 5.10, p. 195}).
This proves that $\imath_\oplus H^1_{\Gamma,s}(\om)\subset H^C_\oplus$.
Assume now that $\Gamma=\Gamma_L$. We shall show that in this case
$v_1\in H^1_0(\omega_1)$. Let us define the function 
$\tilde u\co \om_7\to\R$ by
$$
\tilde u(x,y):=\cases u(x,y)&\text{if $x\in\om_4$}\\0&\text{otherwise}
\endcases
$$
Since $\null^\tau u(x,y)=0$ ${\Cal H}^2$-a.e. on $\Gamma_L$, 
it follows that $\tilde u\in H^1_s(\om_7)$ (here ${\Cal H}^2$ is the
two-dimensional Hausdorff measure in $\R^3$ and $\null^\tau u$ is the
trace of $u$ on $\partial\om$). By Proposition \rft co..,
there exist a null
set $S\subset\R^3$ and a function $v_7\in
H^1(\R^2)$ such that
$$ \alignedat2
\tilde u(x,y)&=v_7(x)&\quad&\text{for all $(x,y)\in\Omega_7\setminus
S$.}\\  \endalignedat
$$
Observe that $v_7=0$ a.e. in $\R^2\setminus\omega_1$.
On the other hand, again by Proposition \rft co..,  we can find a function 
$v_4\in H^1(\omega_1)$ and we can choose the set $S$ in such a way 
that
$$ \alignedat2
u(x,y)&=v_4(x)&\quad&\text{for all $(x,y)\in\Omega_4\setminus
S$.}\\  \endalignedat
$$
It follows that $v_1=v_4=v_7$ almost everywhere in $\omega_1$. This in
turn implies that $\null^\tau v_1(x)=0$
${\Cal H}^1$-a.e. on $\partial\omega_1$ (see again \cite{\rfa alt..}),
i.e. $v_1\in H^1_0(\omega_1)$.
So far, we have proved that $\imath_\oplus (H^1_{\Gamma,s}(\om))\subset
H^C_\oplus$ and, if $\Gamma=\Gamma_L$, $\imath_\oplus
(H^1_{\Gamma,s}(\om))\subset
H^0_\oplus$. 

Assume now that $[v]\in H^C_\oplus$. We shall prove that there exists a
function $u\in H^1_{\Gamma,s}(\om)$, with $\Gamma=\emptyset$, such that
$\imath_\oplus u=[v]$. Let us define a function $u$ on $\om$ in the
following way:
$$
u(x,y):=\cases v_k(x)&\text{if $(x,y)\in\om_k$, $k=1,2,3$}\\
0&\text{otherwise}\endcases
$$
Obviously, $u|_{\om_1}\in H^1(\om_1)$. Moreover, $u|_{\om_5}\in
H^1(\om_5)$. In fact, the function $\tilde v_5\co\omega\to\R$ defined by 
$$
\tilde v_5(x):=\cases
v_1(x)&\text{if $x\in\omega_1$}\\v_2(x)&\text{if $x\in\omega_2$}\\
0&\text{otherwise}
\endcases
$$
is in $H^1(\omega)$, since $\null^\tau v_1(x)=\null^\tau v_2(x)$
${\Cal H}^1$-a.e. on $\partial\omega_2$ (see again \cite{\rfa alt..}). Analogously,
$u|_{\om_6}\in
H^1(\om_6)$. Now since $(\om_l)_{l=1,5,6}$ is an open covering of $\om$,
it follows immediately that $u\in H^1(\om)$. It is easily
verified that $\nabla_y u=0$ almost everywhere, so $u\in
H^1_{\Gamma,s}(\om)$. By construction, $\imath_\oplus u=[v]$. 

Assume now that $[v]\in H^0_\oplus$. We shall prove that there exists a
function $u\in H^1_{\Gamma,s}(\om)$, with $\Gamma=\Gamma_L$, such that
$\imath_\oplus u=[v]$. As before, let us define a function $u$ on $\om$ in
the
following way:
$$
u(x,y):=\cases v_k(x)&\text{if $(x,y)\in\om_k$, $k=1,2,3$}\\
0&\text{otherwise}\endcases
$$
By the same arguments as above, it follows easily that $u\in H^1(\om)$
and that $\nabla_y u=0$ almost everywhere. We shall show that
$\null^\tau u=0$ on $\Gamma_L$. 
To this end, let us choose sequences
$(v_k^n)_{n\in\N}$, $v_k^n\in C^1_0(\omega_k)$, $v_k^n\to v_k$ in
$H^1(\omega_k)$ as $n\to\infty$, $k=1,2,3$, and let us define
$$
u^n(x,y):=\cases v_k^n(x)&\text{if $(x,y)\in\om_k$, $k=1,2,3$}\\
0&\text{otherwise}\endcases
$$
for $n\in\N$. Then $u^n\in C^1(\cl\om)$ and $u^n(x)=0$ on $\Gamma_L$ for
all
$n\in\N$. Moreover, it is easy to verify that $u_n\to u$ in $H^1(\om)$, so
we deduce that $u\in H^1_{\Gamma,s}(\om)$. By construction we have that 
$\imath_\oplus u=[v]$. 
This concludes the proof of
part (2).

Now we prove (1). Let $[v]\in L_\oplus$. We shall prove that there exists
$v\in L^2_{\Gamma,s}(\om)$ such that $\imath_\oplus u=[v]$. Again, we
define a
function $u$ on $\om$ in the following way:
$$
u(x,y):=\cases v_k(x)&\text{if $(x,y)\in\om_k$, $k=1,2,3$}\\
0&\text{otherwise}\endcases
$$
Then $u\in L^2(\om)$. We claim that $u\in L^2_{\Gamma,s}(\om)$, both with
$\Gamma=\Gamma_L$ and with $\Gamma=\emptyset$. This means that $u$ can be
approximated in the $L^2$-norm by functions of $H^1_{\Gamma,s}(\om)$. 
To this end, let us choose sequences
$(v_k^n)_{n\in\N}$, $v_k^n\in C^1_0(\omega_k)$, $v_k^n\to v_k$ in
$L^2(\omega_k)$ as $n\to\infty$, $k=1,2,3$, and let us define
$$
u^n(x,y):=\cases v_k^n(x)&\text{if $(x,y)\in\om_k$, $k=1,2,3$}\\
0&\text{otherwise}\endcases
$$
for $n\in\N$. Then, as in the proof of part (1), $u^n\in
H^1_{\Gamma,s}(\om)$ for
all
$n\in\N$, both with $\Gamma=\Gamma_L$ and with $\Gamma=\emptyset$.
Moreover, it is easy to verify that $u_n\to u$ in $L^2(\om)$, so
we deduce that $u\in L^2_{\Gamma,s}(\om)$. By construction we have that 
$\imath_\oplus u=[v]$ and the proof is complete.
\qed\enddemo

\proclaim{Corollary \dft infinito..}
The space $H^1_{\Gamma,s}(\om)$ is infinite dimensional, both with
$\Gamma=\emptyset$ and with $\Gamma=\Gamma_L$.\qed
\endproclaim

\newsection\head \the\secnumber.
$H^2$-regularity and characterization of $D(A_0)$\endhead

Let us define the bilinear forms
$$
a_k(u,v):=\int_{\omega_k}h_k\nabla u(x)\cdot\nabla v(x)\d x,\quad
u,v\in H_k
$$
on $H_k\times H_k$, $k=1,2,3$, and the bilinear form
$$
a_\oplus([u],[v]):=a_1(u_1,v_1)+a_2(u_2,v_2)+a_3(u_3,v_3),\quad[u],[v]\in
H_\oplus
$$
on $H_\oplus\times H_\oplus$.
Let us indicate by $a^C_\oplus$ and $a^0_\oplus$ the restrictions of
$a_\oplus$ to $H^C_\oplus\times H^C_\oplus$ and $H^0_\oplus\times
H^0_\oplus$ respectively. Let $A^C_\oplus$ (resp. $A^0_\oplus$) be the
self-adjoint operator generated by $a^C_\oplus$ (resp. $a^0_\oplus$) in
$H^C_\oplus$ (resp. $H^0_\oplus$). Finally, let us indicate simply by $a$
the bilinear form $a_0$ on $H^1_{\Gamma,s}(\om)\times H^1_{\Gamma,s}(\om)$
defined in \rff lbil.., and by $A$ the corresponding
self-adjoint
operator $A_0$ defined in \rff lop...
Observe that 
$$
a(u,v)=a_\oplus(\imath_\oplus
u,\imath_\oplus v)
\quad\text{for $u$ and $v\in H^1_{\Gamma,s}(\om)$.}
$$

Assume that $\Gamma=\emptyset$. If $u\in D(A)$, then, for all $v\in
H^1_{\Gamma,s}(\om)$, we have
$$
\langle Au,v\rangle_{L^2(\om)}=a(u,v)=a_\oplus^C(\imath_\oplus
u,\imath_\oplus v).
$$
On the other hand,
$$
\langle Au,v\rangle_{L^2(\om)}=\langle\imath_\oplus
Au,\imath_\oplus v\rangle_{L_\oplus}.
$$
It follows that 
$$
a_\oplus^C(\imath_\oplus
u,\imath_\oplus v)=\langle\imath_\oplus
Au,\imath_\oplus v\rangle_{L_\oplus}
$$
for all $v\in H^1_{\Gamma,s}(\om)$, so $\imath_\oplus u\in D(A_\oplus^C)$
and $A_\oplus^C\imath_\oplus u=\imath_\oplus Au$. Similarly, one can prove
that, whenever $[u]\in D(A_\oplus^C)$, then $\imath_\oplus^{-1}[u]\in
D(A)$ and $A\imath_\oplus^{-1}[u]=\imath_\oplus^{-1}A_\oplus^C[u]$. This
means that $\imath_\oplus$ restricts to an isometry of $D(A)$ onto
$D(A_\oplus^C)$ and that $A=\imath_\oplus^{-1}A_\oplus^C\imath_\oplus$.

In the same way we can prove that, if $\Gamma=\Gamma_L$, then
$\imath_\oplus$ restricts to an isometry of $D(A)$ onto $D(A_\oplus^0)$
and
that $A=\imath_\oplus^{-1}A_\oplus^0\imath_\oplus$.

So the problem of characterizing $D(A)$ reduces to the problem of
characterizing $D(A_\oplus^C)$ and $D(A_\oplus^0)$.

We need the following regularity result:

\proclaim{Proposition \dft reg..}
Let $[u]\in H_\oplus$ and $[w]\in L_\oplus$. Assume that one of the
following properties
holds:
\roster
\item $[u]\in H_\oplus^C$ and
$$
a_\oplus([u],[v])=\langle [w],[v]\rangle_{L_\oplus}\quad\text{for all
$[v]\in H_\oplus^C$};
$$
\item $[u]\in H_\oplus^0$ and
$$
a_\oplus([u],[v])=\langle [w],[v]\rangle_{L_\oplus}\quad\text{for all
$[v]\in H_\oplus^0$}.
$$
\endroster
Then $u_k\in H^2(\omega_k)$ for $k=1,2,3$.
\endproclaim
\demo{Proof} See the Appendix.\qed\enddemo

For $k=1,2,3$ let us define the spaces
$$
Z_k:=H^2(\omega_k)\quad\text{and}\quad Z^0_k:=H^2(\omega_k)\cap
H^1_0(\omega_k).
$$
Moreover, let us define the spaces
$$
Z_\oplus:=Z_1\oplus Z_2\oplus Z_3\quad\text{and}\quad Z_\oplus^0:=
Z^0_1\oplus Z^0_2\oplus Z^0_3.
$$
Then we have the following characterization of $D(A_\oplus^C)$ and
$D(A^0_\oplus)$:

\proclaim{Theorem \dft char2..}
The following properties hold
\roster
\item $D(A_\oplus^0)=Z^0_\oplus$ and
$$
A_\oplus^0[u]=(-\Delta u_1,-\Delta u_2, -\Delta u_3)\quad\text{for
$[u]\in Z^0_\oplus$};
$$
\item $D(A_\oplus^C)=Z^C_\oplus$ and
$$
A_\oplus^C[u]=(-\Delta u_1,-\Delta u_2, -\Delta u_3)\quad\text{for
$[u]\in Z^C_\oplus$},
$$
where $Z_\oplus^C$ is the subspace of $Z_\oplus$ consisting of all
$[u]=(u_1,u_2,u_3)$ satisfying
$$
\null^\tau u_1(x)=\null^\tau u_2(x)=\null^\tau u_3(x)\quad\text{${\Cal
H}^1$-a.e. on $\partial\omega_2$,}
$$ 
$$
\partial_{\nu_1} u_1(x)=0\quad\text{${\Cal
H}^1$-a.e. on $\partial\omega$} 
$$
and
$$
h_1\nabla u_1\cdot\nu_1+h_2\nabla u_2\cdot\nu_2+h_3\nabla u_3\cdot\nu_3
=0\quad\text{${\Cal
H}^1$-a.e. on $\partial\omega_2$,}\tag\dff balcon..
$$
where $\nu_k$ is the outward normal vector field on $\partial\omega_k$
for $k=1,2,3$. We call \rff balcon.. the (Kirchoff type) {\rm balance
condition} on $\partial\omega_2$.
\endroster
\endproclaim
\demo{Proof}
First we prove (1). Let $[u]\in D(A^0_\oplus)$. Then by definition there exists 
$[w]\in L_\oplus $ such that
$$
a_\oplus([u],[v])=\langle [w],[v]\rangle_{L_\oplus}\quad\text{for all
$[v]\in H_\oplus^0$}.
$$
Since by Proposition \rft reg.. $u_k\in H^2(\omega_k)\cap H^1_0(\omega_k)$ for $k=1,2,3$, 
we obtain immediately that $[u]\in Z^0_\oplus$. Moreover,
a simple
integration by parts yields
$$
\gather
-\sum_{k=1}^3\int_{\omega_k}h_kv_k(x)\Delta u_k(x)\d x=
\sum_{k=1}^3\int_{\omega_k}h_k\nabla v_k(x)\cdot\nabla u_k(x)\d x\\
=\sum_{k=1}^3\int_{\omega_k}h_kv_k(x)w_k(x)\d x
\quad\text{for all
$[v]\in H_\oplus^0$.}
\endgather
$$
Choose $[v]=(v_1,0,0)$, with $v_1\in H^1_0(\omega_1)$ arbitrary. Then by definition $[v]\in
H_\oplus^0$. With this choice, we obtain
$$
-\int_{\omega_1}h_1v_1(x)\Delta u_1(x)\d x
=\int_{\omega_1}h_1v_1(x)w_1(x)\d x
\quad\text{for all
$v_1\in H^1_0(\omega_1)$.}
$$  
This implies that $w_1=-\Delta u_1$. In the same way, we obtain that $w_k=-\Delta u_k$ for
$k=1,2,3$, i.e. $A_\oplus^0[u]=(-\Delta u_1,-\Delta u_2, -\Delta u_3)$. 

Assume conversely that 
$[u]\in Z_\oplus^0$. Then integration by parts implies that  
$$
\sum_{k=1}^3\int_{\omega_k}h_k\nabla v_k(x)\cdot\nabla u_k(x)\d x=
-\sum_{k=1}^3\int_{\omega_k}h_kv_k(x)\Delta u_k(x)\d x
\quad\text{for all
$[v]\in H_\oplus^0$.}
$$
Since $(-\Delta u_1,-\Delta u_2, -\Delta u_3)\in L_\oplus$, it follows that
$[u]\in D(A_\oplus^0)$ and the proof of part (1) is complete.

Part (2) is a little more involved. 
 Let $[u]\in D(A^C_\oplus)$. Then by definition there exists 
$[w]\in L_\oplus $ such that
$$
a_\oplus([u],[v])=\langle [w],[v]\rangle_{L_\oplus}\quad\text{for all
$[v]\in H_\oplus^C$}.
$$
By Proposition \rft reg.., $u_k\in H^2(\omega_k)$ for $k=1,2,3$, so 
we obtain immediately that $[u]\in Z_\oplus$. Moreover, since $[u]\in H^C_\oplus$, 
we have of course $\null^\tau u_1(x)=\null^\tau u_2(x)=\null^\tau u_3(x)$
${\Cal H}^1$-a.e. on $\partial\omega_2$. 
A simple
integration by parts yields
$$
\gather
-\sum_{k=1}^3\int_{\omega_k}h_kv_k(x)\Delta u_k(x)\d x
+\sum_{k=1}^3\int_{\partial\omega_k}h_k 
 v_k(x)\nabla u_k(x)\cdot\nu_k(x)\d{\Cal H}^1 x\\
=\sum_{k=1}^3\int_{\omega_k}h_k\nabla v_k(x)\cdot\nabla u_k(x)\d x
=\sum_{k=1}^3\int_{\omega_k}h_kv_k(x)w_k(x)\d x\\
\quad\text{for all
$[v]\in H_\oplus^C$.}\endgather
$$
Choose $[v]=(v_1,0,0)$, with $v_1\in H^1_0(\omega_1)$ arbitrary. Then $[v]\in
H_\oplus^C$. With this choice, we obtain
$$
-\int_{\omega_1}h_1v_1(x)\Delta u_1(x)\d x
=\int_{\omega_1}h_1v_1(x)w_1(x)\d x
\quad\text{for all
$v_1\in H^1_0(\omega_1)$.}
$$  
Since $H^1_0(\omega_1)$ is dense in $L^2(\omega_1)$, we obtain that $w_1=-\Delta u_1$. In the same
way, we obtain that
$w_k=-\Delta u_k$ for
$k=1,2,3$, i.e. $A_\oplus^C[u]=(-\Delta u_1,-\Delta u_2, -\Delta u_3)$.
Now choose $[v]=(v_1,0,0)$ with $\null^\tau v_1=0$ ${\Cal
H}^1$-a.e. on $\partial\omega_2$. Then $[v]\in H_\oplus^C$ and we obtain
$$
\gather
-\int_{\omega_1}h_1v_1(x)\Delta u_1(x)\d x
+\int_{\partial\omega}h_1 
v_1(x)\nabla u_1(x)\cdot\nu_1(x)\d{\Cal H}^1 x\\
=\int_{\omega_1}h_1\nabla v_1(x)\cdot\nabla u_1(x)\d x=
a_\oplus^C([v],[u])
=\langle[v], A^C_\oplus[u]\rangle_{L_\oplus}\\
=-\int_{\omega_1}h_1 v_1(x)\Delta u_1(x)\d x.
\endgather
$$
It follows that
$$
\int_{\partial\omega}h_1 
 v_1(x)\nabla u_1(x)\cdot\nu_1(x)\d{\Cal H}^1 x=0.
$$
Since $\null^\tau v_1$ can be chosen arbitrarily in a dense subspace of
$L^2(\partial\omega)$, we obtain that
$\partial_{\nu_1} u_1(x)=0$ ${\Cal H}^1$-a.e. on $\partial\omega$. 
Finally, choose
$[v]$ in such a way that $\null^\tau v_1(x)=0$ ${\Cal
H}^1$-a.e. on $\partial\omega$. Then we have
$$
\gather
-\sum_{k=1}^3\int_{\omega_k}h_kv_k(x)\Delta u_k(x)\d x
+\sum_{k=1}^3\int_{\partial\omega_2}h_k 
 v_k(x)\nabla u_k(x)\cdot\nu_k(x)\d{\Cal H}^1 x\\
=\sum_{k=1}^3\int_{\omega_k}h_k\nabla v_k(x)\cdot\nabla u_k(x)\d x=
a_\oplus^C([v],[u])
=\langle[v], A^C_\oplus[u]\rangle_{L_\oplus}\\
=-\sum_{k=1}^3\int_{\omega_k}h_kv_k(x)\Delta u_k(x)\d x.
\endgather
$$
It follows that
$$
\sum_{k=1}^3\int_{\partial\omega_2}h_k 
 v_k(x)\nabla u_k(x)\cdot\nu_k(x)\d{\Cal H}^1 x=0.
$$
Since $[v]\in H_\oplus^C$, we have $\null^\tau v_1(x)=\null^\tau v_2(x)=\null^\tau v_3(x)$
${\Cal H}^1$-a.e. on $\partial\omega_2$. Finally, since 
$\null^\tau v_1$ can be chosen arbitrarily in a dense subspace of
$L^2(\partial\omega_2)$, we obtain that
$h_1\nabla u_1\cdot\nu_1+h_2\nabla u_2\cdot\nu_2+h_3\nabla u_3\cdot\nu_3
=0$ ${\Cal H}^1$-a.e. on $\partial\omega_2$, and hence $[u]\in Z_\oplus^C$.

Assume conversely 
that 
$[u]\in Z_\oplus^C$. Then integration by parts implies that  
$$
\gather
\sum_{k=1}^3\int_{\omega_k}h_k\nabla v_k(x)\cdot\nabla u_k(x)\d x\\
=-\sum_{k=1}^3\int_{\omega_k}h_kv_k(x)\Delta u_k(x)\d x
+\sum_{k=1}^3\int_{\partial\omega_k}h_k 
 v_k(x)\nabla u_k(x)\cdot\nu_k(x)\d{\Cal H}^1 x
\\\quad\text{for all
$[v]\in H_\oplus^C$.}
\endgather
$$
Since $[v]\in H_\oplus^C$, we have $\null^\tau v_1(x)=\null^\tau v_2(x)=\null^\tau v_3(x)$
${\Cal
H}^1$-a.e. on $\partial\omega_2$.
Moreover, since $[u]\in Z_\oplus^C$, we have 
$\partial_{\nu_1} u_1(x)=0$ ${\Cal H}^1$-a.e. on $\partial\omega$ and 
$h_1\nabla u_1\cdot\nu_1+h_2\nabla u_2\cdot\nu_2+h_3\nabla u_3\cdot\nu_3
=0$ ${\Cal H}^1$-a.e. on $\partial\omega_2$.
This implies immediately that
$$
\sum_{k=1}^3\int_{\partial\omega_k}h_k 
 v_k(x)\nabla u_k(x)\cdot\nu_k(x)\d {\Cal H}^1 x=0.
$$
Since $(-\Delta u_1,-\Delta u_2, -\Delta u_3)\in L_\oplus$, it follows that
$[u]\in D(A_\oplus^0)$ and the proof is complete.
\qed\enddemo

\remark{Remark}
Thanks to Theorem \rft char2.., one can easily prove that the semiflow
generated by equation \rff e5.. in $H^1_{\Gamma,s}(\om)$ with $\Gamma=\emptyset$ (resp. 
$\Gamma=\Gamma_L$) and the
semiflow generated by equation \rff e6.. (resp. \rff e7..) are conjugate by mean of 
the isometry $\imath_\oplus$.\endremark

\newsection\head \the\secnumber.
An application: computation of the eigenvalues\endhead

In this section we shall explain how the characterization of $A_0$ and of its domain,
obtained in Section 3, can be exploited, in some specific situations, to compute the
eigenvalues of $A_0$. We shall consider the domain $\Omega$ described in Figure 1:
we choose two real numbers $r$ and $R$, $0<r<R$, and we define
$$
\omega:=\{\,x\in\R^2\mid 0\leq |x|^2<R^2\,\},\quad
\omega_2=\omega_3:=\{\,x\in\R^2\mid 0\leq |x|^2<r^2\,\}.
$$ 

First, we observe that, thanks to Theorem \rft char2..,  in the case
$\Gamma=\Gamma_L$ the abstract
eigenvalue problem
$$
A_0u=\lambda u
$$
is equivalent to the system
$$
\cases
-\Delta u_j=\lambda u_j,&x\in\omega_j,~j=1,2,3\\
u_j=0,&x\in\partial\omega_j,~j=1,2,3
\endcases\tag\dff eigp1..
$$
The equations in this system are
completely decoupled, so in this case the
sequence of
the eigenvalues of $A_0$ is just the union of the sequences of eigenvalues
of the three Dirichlet problems considered separately.
These problems can be easily treated in the standard way
by writing the equations in polar coordinates and then using separation
of variables. This is a classical result and we don't discuss it here.

The case $\Gamma=\emptyset$ is more interesting. 
Thanks to Theorem \rft char2.., the abstract
eigenvalue problem
$$
A_0u=\lambda u
$$
is equivalent to the system
$$
\cases
-\Delta u_j=\lambda u_j,&x\in\omega_j,~j=1,2,3\\
u_1(x)=u_2(x)=u_3(x),&|x|=r,\\  
 \partial _{\nu_1} u_1=0,&
|x|=R,\\ 
\sum_{j=1}^3h_j\nabla u_j\cdot\nu_j=0,& |x|=r.
\endcases \tag\dff eigp2..
$$
Also in this case the
computation exploits polar coordinates and separation of variables, but
we have to be a little careful because of the coupling at the 
`interface' $\{\,|x|=r\,\}$. 
Let us write for simplicity $A:=A^C_\oplus$ and let us indicate by $A^\C$
the complexification of $A$. Then $A^\C$ is a self-adjoint operator in the
complex Hilbert space $L^\C_\oplus:=L_\oplus+\i L_\oplus$ with domain
$D(A^\C)=D(A)+\i D(A)$. The action of $A^\C$ is defined in the
obvious way by $A^\C([u]+\i[v]):=Au+\i Av$. The operators $A$ and $A^\C$
have the same eigenvalues with the same multiplicity.
Let $\Phi\co
]0,+\infty[\times]0,2\pi[\to\R^2$,
$\Phi(\rho,\theta)\mapsto(\rho\cos\theta,\rho\sin\theta)$ be the system of
polar coordinates on $\R^2\setminus(\R_+\times\{\,0\,\})$.

Set $I_1:=]r,R[$ and $I_j:=]0,r[$ for $j=2,3$.
We look for eigenvalue-eigenvector pairs $(\lambda, [u])$, where 
$\lambda\geq0$
and $[u]$ has
the form
$$
\gather
[u]=(u_1,u_2,u_3)\quad
\text{with}\\
(u_j\circ\Phi)(\rho,\theta)=v_j(\rho)e^{\i
n\theta},\quad(\rho,\theta)\in I_j\times]0,2\pi[,\quad j=1,2,3.
\tag\dff polfor..\\
\endgather
$$ 
Here $n\in\Z$ and $v_j\co I_j\to\R$ for $j=1,2,3$. 
Let us recall that the Laplacian in two-dimensional polar coordinates
assumes the form
$$
{{\partial^2}\over{\partial\rho^2}}+
{{1}\over{\rho}}{{\partial}\over{\partial\rho}}+
{{1}\over{\rho^2}}{{\partial^2}\over{\partial\theta^2}}.
$$
Let us fix $n\in\Z$. Then an eigenvalue-eigenvector pair of the form \rff
polfor.. must satisfy
$$
\cases
-\left(v_j''+{{1}\over{\rho}}v_j'-{{n^2}\over{\rho^2}}v_j\right)=\lambda
v_j,&\rho\in I_j,~j=1,2,3,\\
\text{$v_2$ and $v_3$ regular at $0$,}&\\
v_1(r)=v_2(r)=v_3(r),&\\  
 v_1'(R)=0,&\\ 
h_1v_1'(r)=h_2v_2'(r)+h_3v_3'(r).& 
\endcases \tag\dff bess1..
$$
If $\lambda=0$ and $n=0$, the space of solutions of \rff bess1.. is
one-dimensional, and is generated by $(v_1,v_2,v_3)=(1,1,1)$. In fact 
a fundamental system of solutions for the equation
$$
v_j''+{{1}\over{\rho}}v_j'=0
$$
is given by
$1$ and $\log \rho$.
If $\lambda=0$ and $n\not=0$, then \rff bess1.. has no non-trivial
solutions. In fact, a fundamental system of solutions for the equation
$$
v_j''+{{1}\over{\rho}}v_j'-{{n^2}\over{\rho^2}}v_j=0
$$
is given by
$\rho^n$ and $\rho^{-n}$.

Assume now that $\lambda\not=0$. Setting
$\tilde v_j(\xi):=v_j(\xi/\sqrt\lambda)$, we transform the equations
$$
-\left(v_j''+{{1}\over{\rho}}v_j'-{{n^2}\over{\rho^2}}v_j\right)=\lambda
v_j,~j=1,2,3\tag\dff bess2..
$$ 
to
$$
\tilde v_j''+{{1}\over{\xi}}\tilde
v_j'+\left(1-{{n^2}\over{\xi^2}}\right)\tilde v_j =0,~j=1,2,3.\tag\dff
bess3..
$$
The latter are Bessel equations of order $|n|$ and, for $j=1,2,3$, a
corresponding
fundamental system of solutions is given by
$J_{|n|}(\xi)$ and $ Y_{|n|}(\xi)$,
where $J_{|n|}$ and $Y_{|n|}$ are the first and the second Bessel function
of
order $|n|$  (see e.g. \cite{\rfa wal..}). 
It follows that a fundamental system of solutions for the equations \rff
bess2.. for $j=1,2,3$ is given by
$$
J_{|n|}(\sqrt\lambda\rho),\quad Y_{|n|}(\sqrt\lambda\rho).
$$
It is well known that $Y_{|n|}$ is singular at $0$. It follows that, for a
given positive $\lambda$, \rff bess1.. admits nontrivial solutions if and
only if we can find real constants $c_i$, $i=1$, \dots, $4$, with 
$(c_1,c_2,c_3,c_4)\not=(0,0,0,0)$, such that 
$$
\cases
c_1 J_{|n|}'(\sqrt\lambda R)+c_4 Y_{|n|}'(\sqrt\lambda R)=0&\\
c_1 J_{|n|}(\sqrt\lambda r)+c_4 Y_{|n|}(\sqrt\lambda
r)=c_2J_{|n|}(\sqrt\lambda r)&\\
c_2J_{|n|}(\sqrt\lambda r)=c_3J_{|n|}(\sqrt\lambda r)&\\
c_1 h_1J_{|n|}'(\sqrt\lambda r)+c_4 h_1Y_{|n|}'(\sqrt\lambda r)=
c_2h_2J_{|n|}'(\sqrt\lambda r)+c_3h_3J_{|n|}'(\sqrt\lambda r).&\\
\endcases\tag\dff siscom..
$$
This is possible if and only if
$$
\det M(n,\lambda,r,R)=0,
$$
where
$$
M(n,\lambda,r,R)=\left(
\matrix
J_{|n|}'(\sqrt\lambda R)&0&0& Y_{|n|}'(\sqrt\lambda R)\\
J_{|n|}(\sqrt\lambda r)&-J_{|n|}(\sqrt\lambda r)&0&Y_{|n|}(\sqrt\lambda
r)\\
0&J_{|n|}(\sqrt\lambda r)&-J_{|n|}(\sqrt\lambda r)&0\\
h_1J_{|n|}'(\sqrt\lambda r)&-h_2J_{|n|}'(\sqrt\lambda
r)&-h_3J_{|n|}'(\sqrt\lambda
r)&h_1Y_{|n|}'(\sqrt\lambda r)\\
\endmatrix\right).
$$
Observe that $\det M(n,\lambda,r,R)$ is an analytic function of
$\lambda>0$. It follows that, for every $n\in\Z$, the zeroes of
$\det M(n,\lambda,r,R)$ in $\R_+$ form a sequence
$$
\lambda_{nm},\quad m=1,2,3,\dots
$$ 
of eigenvalues of $A^\C$ and hence of $A$. Thus we obtain that the set
$$
\{\,\lambda_{nm}|n\in\Z,m=1,2,3,\dots\,\}\cup\{\,0\,\}
$$
is contained in the sequence of the eigenvalues of $A^\C$ and hence of
$A$. The corresponding eigenfunctions can be computed by solving
the system \rff siscom.. with $\lambda=\lambda_{nm}$. If
$(c_1,c_2,c_3,c_4)$ is a nontrivial real solution of \rff siscom.., then
$$
((c_1 J_{|n|}(\sqrt{\lambda_{nm}}\rho)+c_4
Y_{|n|}(\sqrt{\lambda_{nm}}\rho))e^{\i
n\theta},c_2 J_{|n|}(\sqrt{\lambda_{nm}}\rho)e^{\i n\theta},c_3
J_{|n|}(\sqrt{\lambda_{nm}}\rho)e^{\i n\theta})
$$
is an eigenfunction for the eigenvalue $\lambda_{nm}$, expressed in polar
coordinates. Thus, for $n\in\Z$ and $m=1,2,3,\dots$ fixed, we obtain
a finite set of orthonormal eigenfunctions 
$$
\,\{[u]_{nm}^\ell\mid \ell=1,\dots,p(n,m)\,\}
$$
for the eigenvalue $\lambda_{nm}$. Notice that $p(n,m)\leq4$. However, the
multiplicity of $\lambda_{nm}$ can be larger than $p(n,m)$, since we can
have $\lambda_{\bar n\bar m}=\lambda_{nm}$ for some $\bar n\not=n$.

Finally, we claim that all eigenvalues and eigenfunctions of $A^\C$ can be
obtained in this way. 
To this end, for $n\in\Z$ let us first define the space
$$
\gather
(L^\C_\oplus)_n:=\\
\{\,[u]\in L^\C_\oplus\mid
(u_j\circ\Phi)(\rho,\theta)=v_j(\rho)e^{\i n\theta},~v_j\co I_j\to\C,~
(\rho,\theta)\in I_j\times]0,2\pi[,~
j=1,2,3\,\}.\endgather
$$
Observe that a triple of functions $(u_1,u_2,u_3)$, $u_j\co\omega_j\to\C$,
$j=1,2,3$, satisfying $(u_j\circ\Phi)(\rho,\theta)=v_j(\rho)e^{\i
n\theta}$ for some $v_j\co I_j\to\C$,  $(\rho,\theta)\in
I_j\times]0,2\pi[$, $j=1,2,3$, belongs to 
$(L^\C_\oplus)_n$ if and only if 
$$
\int_{I_j}\rho |v_j(\rho)|^2\d\rho<\infty\quad\text{for $j=1,2,3$}.
$$ 
In fact, $\rho=J\Phi(\rho,\theta)$ for
$(\rho,\theta)\in]0,+\infty[\times]0,2\pi[$.
It is also easy to check that
$$
(L^\C_\oplus)_n\perp (L^\C_\oplus)_{\bar n}\quad\text{for $n\not=\bar n$.}
$$
Moreover,
$$
\cl{\bigoplus_{n\in\Z}(L^\C_\oplus)_n}=L^\C_\oplus.
$$
This is true, since $\{\,e^{\i n\theta}\mid n\in\Z\,\}$ is a complete
orthonormal system in $L^2(]0,2\pi[,\C)$.

Write
$$
[u]_{00}:=(\sum_{j=1}^3h_j|\omega_j|)^{-1/2}(1,1,1).
$$
If we show that, for a fixed $n\in\Z$, $n\not=0$, the set
$$
\{\,[u]_{nm}^\ell\mid \ell=1,\dots,p(n,m),~m=1,2,3,\dots\,\}
$$ 
is a complete orthonormal system in $(L^\C_\oplus)_n$ and that the set 
$$
\{\,[u]_{0m}^\ell\mid \ell=1,\dots,p(0,m),~m=1,2,3,\dots\,\}
\cup\{\,[u]_{00}\,\}
$$ 
is a complete orthonormal system in $(L^\C_\oplus)_0$, 
we are done.

Let us define the Hilbert space
$$
\L_\oplus:=\{\,[v]=(v_1,v_2,v_3)\mid \rho^{1/2}v_j(\rho)\in
L^2(I_j,\R),~j=1,2,3\,\}
$$
equipped with the scalar product
$$
\{[v],[\nu]\}_\oplus:=
\sum_{j=1}^3\int_{I_j}h_j\rho v_j(\rho)\nu_j(\rho)\d\rho,\quad
[v],[\nu]\in\L_\oplus.
$$
Set $\L^\C_\oplus=\L_\oplus+\i\L_\oplus$, i.e.
$$
\L_\oplus^\C:=\{\,[v]=(v_1,v_2,v_3)\mid \rho^{1/2}v_j(\rho)\in
L^2(I_j,\C),~j=1,2,3\,\}.
$$
 Moreover, let us define
the isometry
$$
\gather
\jmath\co \L_\oplus^\C\to (L^\C_\oplus)_n\\
(v_1,v_2,v_3)\mapsto
(2\pi)^{-1/2}(w_1\circ\Phi^{-1},w_3\circ\Phi^{-1},w_3\circ\Phi^{-1}),
\endgather
$$
where
$$
w_j(\rho,\phi):=v_j(\rho)e^{\i n\theta},\quad
(\rho,\theta)\in I_j\times]0,2\pi[,\quad j=1,2,3.
$$
It is enough to prove that  the sets
$$
{\Cal B}_n:=
\{\,\jmath^{-1}[u]_{nm}^\ell\mid
\ell=1,\dots,p(n,m),~m=1,2,3,\dots\,\},\quad n\in\Z\setminus\{0\},
$$ 
and 
$$
\gather
{\Cal B}_0:=
\{\,\jmath^{-1}[u]_{0m}^\ell\mid
\ell=1,\dots,p(0,m),~m=1,2,3,\dots\,\}
\cup\{\,\jmath^{-1}[u]_{00}\,\}
\endgather
$$ 
are complete orthonormal systems in $\L_\oplus^\C$.  
Actually, since
$$
\gather
\jmath^{-1}[u]_{nm}^\ell=(v_{nm,1}^\ell,v_{nm,2}^\ell,v_{nm,3}^\ell)
\in\L_\oplus\\
\text{for $\ell=1,\dots,p(n,m)$, $m=1,2,3,\dots$ and for all $n\in\Z$,} 
\endgather
$$
it is enough to prove that ${\Cal B}_n$ and ${\Cal
B}_0$
are complete
orthonormal systems in $\L_\oplus$.

Set $\lambda_{nm}^\ell:=\lambda_{nm}$ for $\ell=1$, \dots, $p(n,m)$,
$m=1,2,3,\dots$, $n\in\Z$.
For $n\not=0$, the set
$$
{\Cal E}_n:=\{\,(\lambda_{nm}^\ell,\jmath^{-1}[u]_{nm}^\ell)\mid
\ell=1,\dots,p(n,m),~m=1,2,3,\dots\,\}
$$
is by construction the set of eigenvalue-eigenvector pairs of the system
$$
\cases
-\left(v_j''+{{1}\over{\rho}}v_j'-{{n^2}\over{\rho^2}}v_j\right)=\lambda
v_j,&\rho\in I_j,~j=1,2,3,\\
\text{$v_2$ and $v_3$ regular at $0$,}&\\
v_1(r)=v_2(r)=v_3(r),&\\  
 v_1'(R)=0,&\\ 
h_1v_1'(r)=h_2v_2'(r)+h_3v_3'(r).& 
\endcases\tag\dff krk1.. 
$$
For $n=0$, the set
$$
{\Cal E}_0:=\{\,(\lambda_{0m}^\ell,\jmath^{-1}[u]_{0m}^\ell)\mid
\ell=1,\dots,p(0,m),~m=1,2,3,\dots\,\}\cup
\{\,(0,\jmath^{-1}[u]_{00})\,\}
$$
is by construction the set of eigenvalue-eigenvector pairs of the system
$$
\cases
-\left(v_j''+{{1}\over{\rho}}v_j'\right)=\lambda
v_j,&\rho\in I_j,~j=1,2,3,\\
\text{$v_2$ and $v_3$ regular at $0$,}&\\
v_1(r)=v_2(r)=v_3(r),&\\  
 v_1'(R)=0,&\\ 
h_1v_1'(r)=h_2v_2'(r)+h_3v_3'(r).&\endcases
\tag\dff krk2..  
$$
Let us define the spaces
$$
\gather
\H^0_\oplus:=\\
\{\,[v]\in\L_\oplus\mid v_j\in
H^1_\loc(I_j),\,\rho^{1/2}v_j'(\rho)\in L^2(I_j), j=1,2,3,\,
v_1(r)=v_2(r)=v_3(r)\,\}
\endgather
$$
and, for $n\in\Z\setminus\{0\}$, 
$$
\H^n_\oplus:=\{\,[v]\in\H^0_\oplus\mid \rho^{-1/2}v_j(\rho)
\in L^2(I_j), j=1,2,3\,\},
$$
equipped with the scalar products
$$
\gather
\{\{[v],[\nu]\}\}^0_\oplus:=
\sum_{j=1}^3\int_{I_j}h_j\rho v_j'(\rho)\nu_j'(\rho)\d\rho
+\sum_{j=1}^3\int_{I_j}h_j\rho v_j(\rho)\nu_j(\rho)\d\rho,\\
[v],[\nu]\in\H^0_\oplus\endgather
$$
and
$$
\gather
\{\{[v],[\nu]\}\}^n_\oplus:=\\
\sum_{j=1}^3\int_{I_j}h_j\rho v_j'(\rho)\nu_j'(\rho)\d\rho
+\sum_{j=1}^3\int_{I_j}h_j{{n^2}\over{\rho}} v_j(\rho)\nu_j(\rho)\d\rho
+\sum_{j=1}^3\int_{I_j}h_j\rho v_j(\rho)\nu_j(\rho)\d\rho,\\
[v],[\nu]\in\H^n_\oplus\endgather
$$
respectively.
Then one can show that $\H^0_\oplus$ and $\H^n_\oplus$,
$n\in\Z\setminus\{0\}$,
 are densely and compactly
imbedded in $\L_\oplus$. 

Let us define the bilinear forms
$$
\gather
a\{[v],[\nu]\}^0_\oplus:=
\sum_{j=1}^3\int_{I_j}h_j\rho v_j'(\rho)\nu_j'(\rho)\d\rho,\\
[v],[\nu]\in\H^0_\oplus\endgather
$$
and
$$
\gather
a\{[v],[\nu]\}^n_\oplus:=
\sum_{j=1}^3\int_{I_j}h_j\rho v'_j(\rho)\nu'_j(\rho)\d\rho
+\sum_{j=1}^3\int_{I_j}h_j{{n^2}\over{\rho}} v_j(\rho)\nu_j(\rho)\d\rho
\\
[v],[\nu]\in\H^n_\oplus\endgather
$$
on $\H^0_\oplus$ and $\H^n_\oplus$ respectively.
Then  we have that the set
${\Cal E}_0$
is the complete set of `proper value~-- proper vector' pairs of  
$$
\cases [v]\in \H^0_\oplus&\\
a\{[v],[\nu]\}^0_\oplus=\lambda \{[v],[\nu]\}_\oplus&\text{for all 
$[\nu]\in \H^0_\oplus$}
\endcases\tag\dff trst1..
$$
Analogously, for all $n\in\Z\setminus\{0\}$, the set
${\Cal E}_n$
is the complete set of `proper value~-- proper vector' pairs of  
$$
\cases [v]\in \H^n_\oplus&\\
a\{[v],[\nu]\}^n_\oplus=\lambda \{[v],[\nu]\}_\oplus&\text{for all 
$[\nu]\in \H^n_\oplus$}
\endcases\tag\dff trst2..
$$
Actually, \rff trst1.. (resp. \rff trst2..) can be considered as 
the `weak formulation' of \rff krk2.. (resp. \rff krk1..).

By the abstract theory of proper values for couples of bilinear forms
(see e.g. \cite{\rfa rav..} or \cite{\rfa wein..}),
we finally obtain that ${\Cal B}_0$ and ${\Cal B}_n$, $n\in\Z\setminus\{0\}$,
are complete orthonormal systems in $\L_\oplus$.

\newsection\head \the\secnumber.
Appendix\endhead

In this appendix we give the 

\demo{Proof of Proposition \rft reg..}
Assume first that (2) holds and remind that $u_k\in H^1_0(\omega_k)$ for $k=1,2,3$. Choose 
$[v]=(v_1,0,0)$, with $v_1\in H^1_0(\omega_1)$ arbitrary. Then by definition $[v]\in
H_\oplus^0$. With this choice, we obtain
$$
\int_{\omega_1}h_1\nabla u_1(x)\cdot\nabla v_1(x)\d x
=\int_{\omega_1}h_1w_1(x)v_1(x)\d x
\quad\text{for all
$v_1\in H^1_0(\omega_1)$.}
$$
Since $\partial\omega_1$ is of class $C^2$, 
the classical regularity
results for the Dirichlet problem apply to the present situation and we get
without any further effort that $u_1\in H^2(\omega_1)$. In the same way, 
we obtain that $ u_k\in
H^2(\omega_k)$ for $k=1,2,3$. 

If (1) holds, the situation is much more complicated: we cannot apply 
directly the 
classical regularity results for elliptic equations, because of the
coupling at the `interface' $\partial\omega_2$. We shall use a partition of
unity on $\cl\omega$ in order to isolate the regions where no coupling
occurs: within these regions we can again apply the classical results.
On the other hand, the partition of unity allows us to `localize' the
analysis on the interface. The main difficulty consists in the fact that
we have to handle with the three functions $u_1$, $u_2$ and $u_3$
simultaneously. Fortunately, the compatibility condition \rff compcond..,
in local coordinates, is invariant under `horizontal' translations. 
Then we
shall exploit the well known method of translations due to L. Nirenberg
and obtain at once $H^2$ regularity of $u_k$, $k=1,2,3$. 

We start by carefully choosing an open covering of $\cl\omega$. 
Since $\partial\omega_2$ is of class $C^2$, we can cover it by a finite number
of open sets $U_1$, $U_2$, \dots, $U_m$, in such a way that, for $i=1$, \dots, $m$, 
there exists a $C^2$ diffeomorphism
$\Phi_i\colon\, ]-1,1[\,\times \,]-1,1[\to U_i$ with the property that
\roster
\item
$\partial\omega_2\cap U_i=\Phi_i(\,]-1,1[\,\times\{0\})$;
\item $\omega_2\cap U_i=\Phi_i(\,]-1,1[\,\times\,]-1,0[\,)$;
\item $\omega_1\cap U_i=\Phi_i(\,]-1,1[\,\times\,]0,1[\,)$.
\endroster
Notice that $U_i\cap \partial\omega=\emptyset$. Moreover, we take $U_0\subset\subset\omega_2$
in such a way that $U_0,U_1$, \dots, $U_m$ form an open covering of $\cl\omega_2$. Notice
that
$U_0\cap\cl\omega_1=\emptyset$.
Finally, we take $U_{m+1}\subset\subset\R^2$ in such a way that
$U_{m+1}\cap\cl\omega_2=\emptyset$ and
$U_1$, \dots, $U_m,U_{m+1}$ form an open covering of $\cl\omega_1$. Then $U_0$, \dots,
$U_{m+1}$ is an open covering of $\cl\omega$.

For $i=0$, \dots, $m+1$, let $\theta_i\in C^\infty_0(\R^2)$, with $\supp
\theta_i\subset U_i$, be a partition of unity on $\cl\omega$, i.e.
$\sum_{i=0}^{m+1} \theta_i\equiv 1$ on $\cl\omega$. Let us observe that
$\sum_{i=0}^m \theta_i\equiv 1$ on $\cl\omega_2=\cl\omega_3$ and
$\sum_{i=1}^{m+1} \theta_i\equiv 1$ on $\cl\omega_1$.
Then we have 
$$
u_1=\sum_{i=1}^{m+1} \theta_i u_1, \quad
u_2=\sum_{i=0}^m \theta_i u_2, \quad\text{and}\quad
u_3=\sum_{i=0}^m \theta_i u_3.
$$
So it is sufficient to show that
$$
\alignedat2
u_{i,1}:=&\theta_iu_1\in H^2(\omega_1)&\quad&\text{for $i=1$, \dots, $m+1$}\\
u_{i,2}:=&\theta_iu_2\in H^2(\omega_2)&\quad&\text{for $i=0$, \dots, $m$}\\
u_{i,3}:=&\theta_iu_3\in H^2(\omega_3)&\quad&\text{for $i=0$, \dots, $m$}\\
\endalignedat
$$
Let us observe that $\supp u_{m+1,1}\subset U_{m+1}\cap\cl\omega_1$, $\supp
u_{0,2}$ and $\supp
u_{0,3}\subset U_0$, and $\supp u_{i,j}\subset
U_i\cap\cl\omega_j$ for $i=1$, \dots, $m$ and
$j=1,2,3$. 

We prove first that $u_{0,2}\in H^2(\omega_2)$ and $u_{0,3}\in
H^2(\omega_3)$, the simplest case.
Let $v_2\in H^1_0(\omega_2)$. We have
$$
\gather
\int_{\omega_2}\nabla u_{0,2}(x)\cdot\nabla v_2(x)\d x\\
=\int_{\omega_2}u_2(x)\nabla \theta_0(x)\cdot\nabla v_2(x)\d x+
\int_{\omega_2}\theta_0(x)\nabla u_2(x)\cdot\nabla v_2(x)\d x\\
=\int_{\omega_2}u_2(x)\nabla \theta_0(x)\cdot\nabla v_2(x)\d x+
\int_{\omega_2}\nabla u_2(x)\cdot\nabla(\theta_0v_2)(x)\d x\\
-\int_{\omega_2}v_2(x)\nabla u_2(x)\cdot\nabla\theta_0(x)\d x.
\endgather
$$
Since $(0,\theta_0v_2,0)\in H^C_\oplus$, we have
$$
\int_{\omega_2}\nabla u_2(x)\cdot\nabla(\theta_0v_2)(x)\d x
=\int_{\omega_2}w_2(x)\theta_0(x)v_2(x)\d x.
$$
Moreover, since $u_2\in H^1(\omega_2)$ and $v_2\in H^1_0(\omega_2)$,
$$
\int_{\omega_2}u_2(x)\nabla \theta_0(x)\cdot\nabla v_2(x)\d x=
-\int_{\omega_2}\diverg(u_2\nabla \theta_0)(x)v_2(x)\d x.
$$
Let us write
$$
\tilde w_2:=-\diverg(u_2\nabla\theta_0)+w_2\theta_0-\nabla
u_2\cdot\nabla\theta_0.
$$
Then $\tilde w_2\in L^2(\omega_2)$ and
$$
\int_{\omega_2}\nabla u_{0,2}(x)\cdot\nabla v_2(x)\d x=
\int_{\omega_2}\tilde w_2(x)v_2(x)\d x.
$$
Since $v_2\in H^1_0(\omega_2)$ is arbitrary, 
we obtain that $u_{0,2}\in H^1_0(\omega_2)$ is a weak solution of
$$
-\Delta u=\tilde w_2~\text{on $\omega_2$},\quad u=0~\text{on
$\partial\omega_2$}.
$$
Then by the standard regularity results for the Dirichlet problem we
obtain that $u_{0,2}\in H^2(\omega_2)$. In the same way we can prove that
$u_{0,3}\in H^2(\omega_2)$.

Next, we consider $u_{m+1,1}$. As we have already mentioned, $\supp
u_{m+1,1}\subset U_{m+1}\cap\cl\omega_1=(U_{m+1}\cap \omega_1)\cup\partial\omega$.
This implies that $\null^\tau u_{m+1,1}=0$ on $\partial\omega_2$.
Let $v_1\in H^1(\omega_1)$, $\null^\tau v_1=0$ on $\partial\omega_2$.
Then we have
$$
\gather
\int_{\omega_1}\nabla u_{m+1,1}(x)\cdot\nabla v_1(x)\d x\\
=\int_{\omega_1}u_1(x)\nabla \theta_{m+1}(x)\cdot\nabla v_1(x)\d x+
\int_{\omega_1}\theta_{m+1}(x)\nabla u_1(x)\cdot\nabla v_1(x)\d x\\
=\int_{\omega_1}u_1(x)\nabla \theta_{m+1}(x)\cdot\nabla v_1(x)\d x+
\int_{\omega_1}\nabla u_1(x)\cdot\nabla(\theta_{m+1}v_1)(x)\d x\\
-\int_{\omega_1}v_1(x)\nabla u_1(x)\cdot\nabla\theta_{m+1}(x)\d x.
\endgather
$$
Since $(\theta_{m+1}v_1,0,0)\in H^C_\oplus$, we have
$$
\int_{\omega_1}\nabla u_1(x)\cdot\nabla(\theta_{m+1}v_1)(x)\d x
=\int_{\omega_1}w_1(x)\theta_{m+1}(x)v_1(x)\d x.
$$
Let us write
$$
\tilde w_1:=w_1\theta_{m+1}-\nabla
u_1\cdot\nabla\theta_{m+1}
\quad\text{and}\quad \tilde W_1:=u_1\nabla\theta_{m+1}.
$$
Then $\tilde w_1\in L^2(\omega_1)$ and $\tilde W_1\in H^1(\omega_1,\R^2)$
and we have
$$
\gather
\int_{\omega_1}\nabla u_{m+1,1}(x)\cdot\nabla v_1(x)\d x=
\int_{\omega_1}\tilde w_1(x)v_1(x)\d x+
\int_{\omega_1}\tilde
W_1(x)\cdot\nabla v_1(x)\d x\\
\text{for all $v_1\in H^1(\omega_1)$ with $\null^\tau v_1=0$ on
$\partial\omega_2$}.
\endgather
$$
Then we can apply the classical regularity results
for elliptic equations with mixed boundary conditions (see
e.g.\cite{\rfa tro..}).
Observe that $\partial\omega_1=\partial\omega\cup\partial\omega_2$ and that the Dirichlet
condition is imposed on the whole $\partial\omega_2$, whereas no a-priori condition
is imposed on $\partial\omega$. Since $\partial\omega_2$ and $\partial\omega$ are smooth and
both closed and open in $\partial\omega_1$, all the hypotheses of Theorem 2.24
in \cite{\rfa tro..} are satisfied. 
So we obtain that $u_{m+1,1}\in H^2(\omega_1)$.

Finally, we shall prove that $u_{i,j}\in H^2(U_i\cap\omega_j)$ for
$j=1,2,3$ and $i=1$, \dots, $m$. Let us fix $i=1$, \dots, $m$, and let us take
$(v_1,v_2,v_3)\in H^C_\oplus$ with $\supp v_{j}\subset U_i\cap\cl\omega_j$ for $j=1,2,3$.
Then we have
$$
\gather
\sum_{j=1}^3\int_{U_i\cap\omega_j}h_j\nabla u_{i,j}(x)\cdot\nabla v_j(x)\d x\\
=\sum_{j=1}^3\int_{U_i\cap\omega_j}h_ju_j(x)\nabla \theta_i(x)\cdot\nabla v_j(x)\d x+
\sum_{j=1}^3\int_{U_i\cap\omega_j}h_j\theta_i(x)\nabla u_j(x)\cdot\nabla v_j(x)\d x\\
=\sum_{j=1}^3\int_{U_i\cap\omega_j}h_ju_j(x)\nabla \theta_i(x)\cdot\nabla v_j(x)\d x+
\sum_{j=1}^3\int_{U_i\cap\omega_j}h_j\nabla u_j(x)\cdot\nabla(\theta_iv_j)(x)\d x\\
-\sum_{j=1}^3\int_{U_i\cap\omega_j}h_jv_j(x)\nabla u_j(x)\cdot\nabla\theta_i(x)\d x.
\endgather
$$
Now observe that $(\theta_iv_1,\theta_iv_2,\theta_iv_3)\in H^C_\oplus$, so
$$
\gather
\sum_{j=1}^3\int_{U_i\cap\omega_j}h_j\nabla u_j(x)\cdot\nabla(\theta_iv_j)(x)\d x
=\sum_{j=1}^3\int_{\omega_j}h_j\nabla u_j(x)\cdot\nabla(\theta_iv_j)(x)\d x\\
=\sum_{j=1}^3\int_{\omega_j}h_jw_j(x)\theta_i(x)v_j(x)\d x
=\sum_{j=1}^3\int_{U_i\cap\omega_j}h_jw_j(x)\theta_i(x)v_j(x)\d x\endgather
$$
Let us write
$$
\tilde w_j:=w_j\theta_i-\nabla
u_j\cdot\nabla\theta_i
\quad\text{and}\quad \tilde W_j:=u_j\nabla\theta_i\quad\text{for $j=1,2,3$.}
$$
Then $\tilde w_j\in L^2(\omega_j)$ and $\tilde W_j\in H^1(\omega_j,\R^2)$ for $j=1,2,3$,
and we have
$$
\gather
\sum_{j=1}^3\int_{U_i\cap\omega_j}h_j\nabla u_{i,j}(x)\cdot\nabla v_j(x)\d x\\
=\sum_{j=1}^3\int_{U_i\cap\omega_j}h_j\tilde w_j(x)v_j(x)\d x+
\sum_{j=1}^3\int_{U_i\cap\omega_j}h_j\tilde W_j(x)\cdot\nabla v_j(x)\d x\tag\dff abrtgz..\\
\text{for all $[v]\in H^C_\oplus$ with $\supp v_{j}\subset 
U_i\cap\cl\omega_j$ for $j=1,2,3$.}
\endgather
$$
 
Set $Q_i:=]-1,1[\times]-1,1[$,
$$
\gather
Q_{i,j}:=\Phi_i^{-1}(U_i\cap\omega_j)\quad\text{for $j=1,2,3$,}\\
\text{i.e.} \quad
Q_{i,1}=]-1,1[\times]0,1[,\quad
Q_{i,2}=Q_{i,3}=]-1,1[\times]-1,0[, 
\endgather
$$
and
$$
\gather
\bar u_{i,j}(\xi):=u_{i,j}(\Phi(\xi)),\quad \bar v_j(\xi):=v_j(\Phi_i(\xi))\\
\text{for $\xi\in Q_{i,j}$, $j=1,2,3$.}
\endgather
$$
Then $\bar u_{i,j}$ and $\bar v_j\in H^1(Q_{i,j})$. 
Moreover, $\supp \bar u_{i,j}$ and $\supp\bar
v_j$ are contained in $Q_{i,j}\cup(\,]-1,1[\times\{0\})$.
Besides, 
$\null^\tau \bar u_{i,1}=\null^\tau \bar
u_{i,2}=\null^\tau \bar u_{i,3}$ and 
$\null^\tau \bar v_1=\null^\tau \bar v_2=\null^\tau \bar v_3$
${\Cal H}^1$-almost everywhere on  $]-1,1[\times\{0\}$. 
Then, changing coordinates in \rff abrtgz.., we have
$$
\gather
\sum_{j=1}^3\int_{Q_{i,j}}h_j
J\Phi_i(\xi)D\Phi_i^{-1}(\Phi_i(\xi))D\Phi_i^{-1}(\Phi_i(\xi))^T\nabla
\bar u_{i,j}(\xi)\cdot\nabla\bar v_j(\xi)\d \xi\\ 
=\sum_{j=1}^3\int_{Q_{i,j}}h_jJ\Phi_i(\xi)\bar
w_j(\xi)\bar v_j(\xi)\d \xi\\
+\sum_{j=1}^3\int_{Q_{i,j}}h_jJ\Phi_i(\xi)D\Phi_i^{-1}(\Phi_i(\xi))
\bar W_j(\xi)\cdot\nabla \bar v_j(\xi)\d \xi,\\
\endgather
$$
where $J\Phi_i(\xi)$ is the Jacobian determinant of $D\Phi_i(\xi)$, $\bar w_j(\xi):=\tilde
w_j(\Phi_i(\xi))\in L^2(Q_{i,j})$ and 
$\bar W_j(\xi):=\tilde W_j(\Phi_i(\xi))\in H^1(Q_{i,j},\R^2)$ for $j=1,2,3$.
Write
$$
J\Phi_i D\Phi_i^{-1}(\Phi_i)D\Phi_i^{-1}(\Phi_i)^T=:(g^{\mu\nu}_i)_{\mu\nu}\in C^1(
Q_i,{\Cal M}(2\times 2)),
$$
$$
J\Phi_i\bar w_j=:\alpha_j\in L^2(Q_{i,j})\quad\text{and}\quad
J\Phi_i D\Phi_i^{-1}\bar W_j=:\beta_j\in H^1(Q_{i,j},\R^2).
$$
Observe also that the matrix $(g^{\mu\nu}_i)_{\mu\nu}$ is symmetric and 
uniformly strongly elliptic on
$ Q_i$, i.e. there exists a positive constant $K$ such that
$$
\sum_{\mu,\nu=1}^2g^{\mu\nu}_i(\xi)h_\mu h_\nu\geq K|h|^2\quad\text{for all $\xi\in\cl Q_i$
and all
$h\in\R^2$.}
$$
Then we have
$$
\gather
\sum_{j=1}^3\int_{Q_{i,j}}h_j \sum_{\mu,\nu=1}^2g^{\mu\nu}_i(\xi)
\partial_\mu\bar u_{i,j}(\xi)\partial_\nu\bar v_j(\xi)\d\xi\\  
=\sum_{j=1}^3\int_{Q_{i,j}}h_j\alpha_j(\xi)\bar v_j(\xi)\d \xi+
\sum_{j=1}^3\int_{Q_{i,j}}h_j\beta_j(\xi)\cdot\nabla \bar v_j(\xi)\d \xi
\tag\dff weak..\endgather
$$
for all $(\bar v_1,\bar v_2,\bar v_3)\in H^C_\oplus(Q_i)$, 
where $H^C_\oplus(Q_i)$ is the set of all
triples   $(\bar v_1,\bar v_2,\bar v_3)\in
H^1(Q_{i,1})\times H^1(Q_{i,2})\times H^1(Q_{i,3})$ with
$\supp \bar v_j\subset Q_{i,j}\cup(]-1,1[\times\{0\})$ and  
$\null^\tau \bar v_1=\null^\tau \bar
v_2=\null^\tau \bar v_3$
${\Cal H}^1$-almost everywhere on  $]-1,1[\times\{0\}$.

Now we are in a position to use the method of translations of Nirenberg. 
First, let us recall that for $u\in L^1_\loc(\R^n)$ and
$h\in\R^n$ one defines
$$
\tau_h u(z):=u(z+h)\quad\text{ and}\quad \delta_h u(z):={{\tau_h u(z)-u(z)}\over{h}},
\quad\text{for $z\in\R^n$}.
$$
We shall use `horizontal' translations: let $h:=(\chi,0)\in\R^2$, with  
$$
|h|<(1/2)\dist (\supp \bar u_{i,j},\{\,-1,1\,\}\times\R)\quad\text{for $j=1,2,3$.}
$$
Then it is very easy to see that $(\tau_h\bar u_{i,1}, \tau_h\bar u_{i,2}, \tau_h\bar
u_{i,3})$,
$(\delta_h\bar u_{i,1}, \delta_h\bar u_{i,2}, \delta_h\bar u_{i,3})$ and
$(\delta_{-h}\delta_h\bar u_{i,1}, \delta_{-h}\delta_h\bar u_{i,2}, 
\delta_{-h}\delta_h\bar u_{i,3})\in H^C_\oplus(Q_i)$. So we can use
$$
(\bar v_1,\bar v_2,\bar v_3):=(\delta_{-h}\delta_h\bar u_{i,1}, \delta_{-h}\delta_h\bar
u_{i,2}, 
\delta_{-h}\delta_h\bar u_{i,3})
$$
as a test function in \rff weak... A simple change of variable yields
$$
\gather
\sum_{j=1}^3\int_{Q_{i,j}}h_j \sum_{\mu,\nu=1}^2\delta_h (g^{\mu\nu}_i
\partial_\mu\bar u_{i,j})(\xi)\partial_\nu(\delta_h\bar u_{i,j})(\xi)\d\xi\\  
=-\sum_{j=1}^3\int_{Q_{i,j}}h_j\alpha_j(\xi) \delta_{-h}\delta_h
\bar u_{i,j}(\xi)\d \xi+
\sum_{j=1}^3\int_{Q_{i,j}}h_j(\delta_h\beta_j)(\xi)\cdot
\nabla(\delta_h \bar u_{i,j})(\xi)\d \xi.
\endgather
$$
Since
$$
\delta_h (g^{\mu\nu}_i\partial_\mu\bar u_{i,j})=
\tau_h (g^{\mu\nu}_i)\partial_\mu(\delta_h\bar u_{i,j})+
\delta_h (g^{\mu\nu}_i)\partial_\mu\bar u_{i,j},
$$
we obtain
$$
\gather
\sum_{j=1}^3\int_{Q_{i,j}}h_j \sum_{\mu,\nu=1}^2(\tau_h g^{\mu\nu}_i)(\xi)
\partial_\mu(\delta_h\bar u_{i,j})(\xi)
\partial_\nu(\delta_h\bar u_{i,j})(\xi)\d\xi\\  
=-\sum_{j=1}^3\int_{Q_{i,j}}h_j \sum_{\mu,\nu=1}^2(\delta_h g^{\mu\nu}_i)(\xi)
\partial_\mu\bar u_{i,j}(\xi)\partial_\nu(\delta_h\bar u_{i,j})(\xi)\d\xi\\ 
-\sum_{j=1}^3\int_{Q_{i,j}}h_j\alpha_j(\xi) \delta_{-h}\delta_h\bar
u_{i,j}(\xi)\d \xi+
\sum_{j=1}^3\int_{Q_{i,j}}h_j(\delta_h\beta_j)(\xi)\cdot
\nabla(\delta_h \bar u_{i,j})(\xi)\d \xi.
\endgather
$$
Now let us recall that
$$
|\delta_{-h}\delta_h\bar u_{i,j}|_{L^2(Q_{i,j})}\leq|\nabla(\delta_h \bar
u_{ij})|_{L^2(Q_{i,j},\R^2)}
$$
and
$$
|\delta_h\beta_j|_{L^2(Q_{i,j},\R^2)}\leq|D\beta|_{L^2(Q_{i,j},{\Cal M}(2\times2))}.
$$
So we get
$$
\alignedat1
\bar K\sum_{j=1}^3\int_{Q_{i,j}}|\nabla(\delta_h\bar u_{i,j})|^2\d\xi&\leq
\sum_{j=1}^3|(g^{\mu\nu}_i)|_{C^1(\cl Q_i)}|\nabla\bar u_{i,j}|_{L^2(Q_{i,j},\R^2)}
|\nabla(\delta_h \bar u_{i,j})|_{L^2(Q_{i,j},\R^2)}\\
&+\sum_{j=1}^3|\alpha_j|_{L^2(Q_{i,j})}|\nabla(\delta_h \bar u_{i,j})|_{L^2(Q_{i,j},\R^2)}\\
&+\sum_{j=1}^3|D\beta|_{L^2(Q_{i,j},{\Cal M}(2\times2))}
|\nabla(\delta_h \bar u_{i,j})|_{L^2(Q_{i,j},\R^2)},
\endalignedat
$$
for some positive constant $\bar K$. This in turn implies that there exists a constant $C>0$
such that
$$
\sum_{j=1}^3|\nabla(\delta_h \bar u_{i,j})|_{L^2(Q_{i,j},\R^2)}^2\leq C
\sum_{j=1}^3|\nabla(\delta_h \bar u_{i,j})|_{L^2(Q_{i,j},\R^2)}
$$
and hence
$$
\left(\sum_{j=1}^3|\nabla(\delta_h \bar u_{i,j})|_{L^2(Q_{i,j},\R^2)}\right)^2\leq 3C
\sum_{j=1}^3|\nabla(\delta_h \bar u_{i,j})|_{L^2(Q_{i,j},\R^2)}.
$$
So, for all sufficiently small $h=(\chi,0)$, we have obtained that 
$$
|\delta_h(\nabla \bar u_{i,j})|_{L^2(Q_{i,j},\R^2)}\leq 
3C\quad\text{for $j=1,2,3$.}\tag\dff
stima..
$$
It is well known that estimates \rff stima.. hold if and only if 
$$
\partial_1\partial_\nu \bar u_{i,j}\in L^2(Q_{i,j})\quad\text{for $\nu=1,2$ and for
$j=1,2,3$}
$$
So, in order to complete the proof, we only need to show that 
$\partial^2_2 \bar u_{i,j}\in L^2(Q_{i,j})$ for $j=1,2,3$. 
This can be easily done by mean of straightforward manipulations of the
distributional identities
$$
-\sum_{\mu,\nu=1}^2\partial_\nu(g^{\mu\nu}_i\partial_\mu \bar u_{i,j})=
\alpha_j-\sum_{\nu=1}^2\partial_\nu\beta_j^\nu,\quad j=1,2,3,
$$  
like in the classical proof of regularity for elliptic equations.  
The proof is complete.
\qed\enddemo

\enddocument


\Refs

\ref\no \dfa alt.. \by H. W. Alt\book Lineare
Funktionalanalysis\bookinfo 2. Auflage\publ
Springer-Verlag\publaddr Berlin Heidelberg New York\yr
1992\endref


\ref \no \dfa arri1..\by J. Arrieta\paper Neumann eigenvalue
problems on exterior perturbations of the domain \jour J. Differ.
Equations \vol 118 \yr1995\pages 54 -- 103
\endref



\ref \no \dfa arriha..\by J. Arrieta, J. Hale and Q. Han\paper
Eigenvalue problems for nonsmoothly perturbed domains \jour J.
Differ. Equations \vol 91 \yr1991\pages 24 -- 52
\endref





\ref \no \dfa carry.. \by M. C. Carbinatto and K. P. Rybakowski
\paper Conley index continuation and thin domain problems \jour Topological
Methods in Nonlinear Analysis \vol 16 \yr 2000 \pages 201--251\endref



\ref \no \dfa Ciu1..\by I.S. Ciuperca\paper Spectral properties of
Schr\"odinger operators on domains with varying order of thinness
\jour J. Dyn. Differ. Equations \vol10 \yr1998\pages 73 -- 108
\endref

\ref\no \dfa hala.. \by Jack K. Hale\book Asymptotic Behavior of
Dissipative Systems\bookinfo Math. Surveys Monographs 25\publ
AMS\publaddr Providence\yr 1988\endref

\ref \no \dfa HaRau1..\by J. Hale and G. Raugel\paper
Reaction-diffusion equations on thin domains \jour J. Math. Pures
Appl. \vol IX Ser.71 \yr1992\pages 33 -- 95
\endref

\ref \no \dfa HaRau2..\by J. Hale and G. Raugel\paper A damped
hyperbolic equation on thin domains \jour  Trans. Am. Math. Soc.
\vol 329 \yr1992\pages 185 -- 219
\endref

\ref \no \dfa HaRau3..\by J. Hale and G. Raugel\paper A
reaction-diffusion equation on a thin $L$-shaped domain \jour
Proc. Roy. Soc. Edinb., Sect A \vol 125 \yr1995\pages 283 -- 327
\endref

\ref\no \dfa He.. \by D. Henry \book Geometric Theory of
Semilinear Parabolic Equations \bookinfo Lecture notes in
mathematics, Vol 840 \publ Springer-Verlag \publaddr NY \yr
1981\endref

\ref\no \dfa lady.. \by O. Ladyzhenskaya \book Attractors for
Semigroups and Evolution Equations \publ Cambridge University
Press \publaddr Cambridge \yr 1991\endref

\ref \no \dfa priry..\by M. Prizzi, M. Rinaldi and K. P. Rybakowski\paper 
Curved thin domains and parabolic equations \jour Studia Mathematica \toappear 
\endref

\ref \no \dfa pr..\by M. Prizzi and K. P. Rybakowski\paper The
effect of domain squeezing upon the dynamics of reaction-diffusion
equations\jour Journal of Differential Equations \vol 173 \yr 2001
\pages 271--320
\endref

\ref \no \dfa pr2..\by M. Prizzi and K. P. Rybakowski\paper 
Inertial manifolds on squeezed domains \jour J. Dynam. Differential
Equations \toappear
\endref


\ref \no \dfa pr3..\by M. Prizzi and K. P. Rybakowski\paper 
Some Recent Results on Thin Domain Problems\jour Topological
Methods in Nonlinear Analysis\vol 14\yr 1999\pages 239--255
\endref


\ref\no \dfa Rau1.. \by G. Raugel \book Dynamics of partial
differential equations on thin domains \bookinfo R. Johnson (ed.),
Dynamical systems. Lectures given at the 2nd session of the Centro
Internazionale Matematico Estivo (CIME) held in Montecatini Terme,
Italy, June 13 -- 22, 1994. Lecture notes in mathematics, Vol.
1609 \publ Springer-Verlag \publaddr Berlin \yr 1995 \pages 208 --
315\endref

\ref\no \dfa rav.. \by P. A. Raviart and J. M.
Thomas\book Introduction \`a l'Analyse Num\'erique des \'Equations
aux D\'eriv\'ees Partielles\publ Masson\publaddr Paris\yr
1983\endref


\ref\no \dfa tro.. \by G.M. Troianiello \book Elliptic Differential
Equations and Obstacle Problems \publ Plenum
Press \publaddr New York \yr 1987\endref

\ref\no \dfa wal.. \by W. Walter \book Gew\"ohnliche
Differentialgleichungen\publ Springer Verlag \publaddr Berlin \yr
1993\endref


\ref\no \dfa wein.. \by Hans F. Weinberger\book
Variational Methods for Eigenvalue Approximation\bookinfo CBMS -
NSF Regional Conference Series in Applied Mathematics\publ
SIAM\publaddr Philadelphia\yr 1974\endref

\endRefs
\bye